\documentclass[a4paper,11pt] {article}
\usepackage{amsmath,amsthm}
\usepackage[english]{inputenc}
\usepackage{geometry}
\usepackage{mdwlist}
\newtheorem{theo}{Theorem}[section]
\numberwithin{equation}{theo}

\numberwithin{equation}{propo}
\usepackage{amssymb}
\usepackage{tikz}
\usepackage{titlesec}
\usetikzlibrary{decorations}
\usetikzlibrary{decorations.markings,arrows}
\usepackage{cite}
\usepackage{comment}
\newtheorem{rmq}{Remark}
\newtheorem{df}{Definition}
\newtheorem{thm}{Theorem}
\newtheorem*{thm1*}{Theorem 1}
\newtheorem*{thm2*}{Theorem 2}
\newtheorem*{thm3*}{Theorem 3}
\newtheorem{prop}{Proposition}
\newtheorem*{thm*}{Theorem}
\newtheorem*{prop*}{Proposition}
\newtheorem*{prop2*}{Proposition 2}
\newtheorem*{prop3*}{Proposition 3}
\newtheorem*{prop5*}{Proposition 5}
\newtheorem{lm}{Lemma}
\newtheorem{conj}{Conjecture}

\titleformat\section{}{}{0pt}{\Large\scshape\bfseries\filcenter\thesection{} - }

\numberwithin{equation}{section}

\begin{document}

\title{Complex behaviour in cyclic competition bimatrix games}
\author{Cezary Olszowiec}
\date{}

\maketitle



\noindent
\let\thefootnote\relax\footnote{{
\hspace{-0.6cm}c.olszowiec14@imperial.ac.uk \\
\hspace{-0.6cm}Department of Mathematics, Imperial College London \\
\hspace{-0.6cm}180 Queen's Gate, SW7 2AZ London \\
}}

{\bf Abstract.} We consider an example of cyclic competition bimatrix game which is a Rock-Scissors-Paper game with assumption about perfect memory of the playing agents. At first we investigate the dynamics in the neighbourhood of the Nash equilibrium as well as the dynamics on the boundary of codimension 1 - that is when one of the strategies is not played by any agent. For the analysis of the asymptotic behaviour close to the boundary of state space, we provide the description of a naturally appearing heteroclinic network. Due to the symmetry, the heteroclinic network induces a quotient network. This quotient network is investigated as well. In the literature this model was already studied, but unfortunately the results stated there are invalid, as we will show in this paper. It turns out that certain types of behaviour are never possible or appear in the system only for some parameter values, contradicting what was stated before in the literature. Moreover, the parameter space $(-1,1)^2$ is divided into four regions where we observe either irregular behaviour, or preference to follow an itinerary consisting of strategies for which one or the other agent does not lose, or they alternate in winning. These regions in parameter space are separated by two analytical curves and lines where the game is either symmetric (system is not $C^1$ linearisable at stationary saddle points) or is zero-sum and the system is Hamiltonian. On each of these curves we observe different bifurcation scenarios: e.g. transition from order to chaos, or from one kind of stability to another kind, or just loss of one dimension of the local stable manifold of the subcycle.

\section{Introduction}
In the last years there was a huge interest in bimatrix games, learning games, models with applications to population dynamics and economics. This paper will consider one such model, that is a Rock-Scissors-Paper game described by the replicator equations, which is a simplified model of a cyclic competition bimatrix game. This model was derived in e.g. \cite{[18]}. Here we will mainly refer to four papers on this topic:
Aguiar M.A.D., Castro S.B.S.D "Chaotic switching in a two-person game", (Physica D 239 (2010) 1598-1609)
and three others by Sato et al. \cite{[18]}, \cite{[16]}, \cite{[17]}.
\newline 
We suggest introductions of mentioned papers to the readers interested in the history of this particular problem and bimatrix games in general.
\newline 
In the paper by Aguiar and Castro, the authors considered the same model as we do. 
However, in their paroofs and their analysis, they made a "simplifying assumption" namely that a certain transition maps between sections associated to stationary points in the heteroclinic network are "identity maps" in the natural coordinates. As we will see in this paper this assumption is not valid in this model. It turns out that when uses the correct transition map, the results one obtains are rather different from those stated in Aguiar\&Castro. A simple consideration based on the mean value theorem shows that it is impossible to assume that these transition maps are "identity maps". 
\newline
In proposition 9, we compute the linear part of the transition maps in convenient local coordinates. This together with lemma 2, shows that the transition maps are in fact affine mappings.  In proposition 12 we present the approximation of the composition of local maps with transition maps after rescaling of the coordinates, the final form of these compositions resemble expressions from Aguiar\& Castro paper.  
\newline
In theorem 4.5 of Aguiar\&Castro, they claim that there is infinite switching near the heteroclinic network $\Sigma_{\Gamma}$. However as we will show in this paper, this result is incorrect. The first reason for this is the following: 
\begin{thm3*}
For every $\epsilon_X,\epsilon_Y \in (-1,1)$ the finite itineraries of the form $tie \rightarrow win/loss \rightarrow loss/win \rightarrow tie$ and $tie \rightarrow loss/win \rightarrow win/loss \rightarrow tie$ are not attainable for $h$ - sufficiently small.
\end{thm3*}
Moreover, in section \ref{sledzenie} we will provide examples showing that attainability of some finite itineraries in arbitrarily small neighbourhood of the heteroclinic network strongly depends on the parameter values of the system.
\newline 
Theorem 4.8, from Aguiar\&Castro paper, states that for $\epsilon_X+\epsilon_Y<0$ the $C_0$ cycle is relatively asymptotically stable (and cannot be essentially asymptotically stable), while for $\epsilon_X+\epsilon_Y>0$ cycles $C_1, C_2$ are relatively asymptotically stable. 
We correct this statements in the theorems \ref{C0eas} and \ref{C2attract}:
\begin{thm1*} 
For $\epsilon_X+\epsilon_Y<0$, the cycle $C_0$ is essentially asymptotically stable.
\end{thm1*}
\begin{thm2*} 
For some open regions 
in parameter space $(\epsilon_X, \epsilon_Y)$, cycles $C_1$ and $C_2$, respectively, are neither almost completely unstable nor essentially asymptotically stable. For each of them, there exist open set of points in the phase space which is attracted to the respective cycle.
\end{thm2*}
In the papers by Sato et al. \cite{[18]}, \cite{[16]}, \cite{[17]}, the derivation of the equations for the RSP game is presented. As well - these equations are transformed into other coordinates which are 
numerically stable, hence very useful for the numerical simulations. The authors gave partial and vague description of the global behaviour in the phase space, depending on the parameter values $\epsilon_X$ and $\epsilon_Y$ based on the numerical investigations. However they mainly concentrate on the case $\epsilon_X+\epsilon_Y=0$. The numerical results reported in these papers are not compatible with each other and also with our results.

\subsection{Organisation of the paper}
We divide this paper into 7 sections.
\newline
In the second section we give the introduction to the model, characterize Nash equilibrium and we explain the global behaviour on the boundaries of codimensions $2$ and $1$.
\newline
Following section provides full construction of the heteroclinic network naturally appearing in the system as well as all of the transition maps associated to this network.
\newline
Asymptotic analysis, stability of the cycles within the network and bifurcations happening in the system are described in section 4. We also present the link between our model and the one introduced by Aguiar\&Castro, and correct some of their results.
\newline
Section number five, consists of examples showing some phenomenon in the system - attainability of the itineraries along heteroclinic orbits as a dependence of parameter values $(\epsilon_X$, $\epsilon_Y)$. We present our examples contradicting the main results from Aguiar$\&$Castro paper and start explaining our presumptions about the global behaviour in the system as well as local behaviour near the heteroclinic network.
\newline
In the sixth section we present our numerical investigations and compare it to Sato's et al simulations \cite{[16]}, \cite{[17]}, \cite{[18]}. We present our presumptions about the way the chaos might be born, we explain the obstacles that forbid us from proving our conjectures. 
\newline
Finally, in the last section we make conclusions, discuss open questions and future work.

\section{Preliminaries}\label{s:2}
We consider the Rock-Scissors-Paper game, played by two agents X and Y. 
\newline
Let $x_1, x_2, x_3 \geq 0$ $( \ y_1, y_2, y_3 \geq 0 \ )$ denote the probability that player X (resp. Y) plays Rock, Scissors or Paper. Obviously $x_1 + x_2 + x_3 =1$, $y_1 + y_2 + y_3 =1$. For each player, the state space is a two-dimensional simplex, hence the whole state space for played game is $\Delta = \Delta_X \times \Delta_Y$. The interaction matrices after normalization are given by:
\begin{equation}\label{macierzeAB}
\begin{aligned}
A= \left( \begin{array}{ccc}
\frac{2}{3} \epsilon_X & 1-\frac{1}{3}\epsilon_X & -1-\frac{1}{3}\epsilon_X \\
-1-\frac{1}{3}\epsilon_X & \frac{2}{3}\epsilon_X & 1-\frac{1}{3}\epsilon_X \\
1-\frac{1}{3}\epsilon_X & -1-\frac{1}{3}\epsilon_X & \frac{2}{3}\epsilon_X \end{array} \right) 
\\
B= \left( \begin{array}{ccc}
\frac{2}{3} \epsilon_Y & 1-\frac{1}{3}\epsilon_Y & -1-\frac{1}{3}\epsilon_Y \\
-1-\frac{1}{3}\epsilon_Y & \frac{2}{3}\epsilon_Y & 1-\frac{1}{3}\epsilon_Y \\
1-\frac{1}{3}\epsilon_Y & -1-\frac{1}{3}\epsilon_Y & \frac{2}{3}\epsilon_Y \end{array} \right)
\end{aligned}
\end{equation}
\newline
with $\epsilon_X,\epsilon_Y \in (-1,1)$, the rewards for ties.
\newline 
Assuming perfect memory of players, the dynamics are described by the replicator equations (see \cite{[18]}):
\begin{equation}\label{rownanie}
\begin{cases}
\dot{x_i}=x_i [(Ay)_i-x^T Ay]
\\
\dot{y_j}=y_j[(Bx)_j-y^T Bx]
\end{cases}
\end{equation}
\newline
with $i,j=1,2,3$. 
\newline
After substitutions $x_3=1-x_1-x_2, y_3=1-y_1-y_2$, we will consider it as a differential equation: 
\begin{equation}\label{rownanieZredukowane}
\begin{cases}
\dot{x_1}=x_1 [(Ay)_1-x^T Ay] 
\\ 
\dot{x_2}=x_2 [(Ay)_2-x^T Ay]
\\ 
\dot{y_1}=y_1[(Bx)_1-y^T Bx] 
\\
\dot{y_2}=y_2[(Bx)_2-y^T Bx]
\end{cases}
\end{equation}
with constraints $x_1,x_2 \geq 0, x_1+x_2 \leq 1$, $y_1,y_2 \geq 0, y_1+y_2 \leq 1$.
\newline
Let us denote the canonical basis in $\mathbb{R}^4$ by $e_1=(1,0,0,0), e_2=(0,1,0,0), e_3=(0,0,1,0), e_4=(0,0,0,1)$.

\begin{prop}
The system (\ref{rownanie}) is invariant under action of the group generated by the cycle $\sigma$ of order $3$:
\begin{align}
\sigma(x_1,x_2,x_3,y_1,y_2,y_3):=(x_2,x_3,x_1,y_2,y_3,y_1)
\end{align}
\end{prop}

Equation (\ref{rownanieZredukowane}) possesses following equilibria:
\begin{prop2*} 
Nash equilibrium of the system  $(x^*,y^*)=(\frac{1}{3},\frac{1}{3},\frac{1}{3},\frac{1}{3})$
\end{prop2*}
\begin{prop3*}
6 equilibria $Z^a,Z^b,Z^c,Z^d,Z^e,Z^f$ of (1,2,1) type (1-dim. stable manifold, 2-dim. center manifold, 1-dim. unstable manifold), which are centers within the 2 dimensional invariant subspaces.
\end{prop3*}
\begin{prop5*}
9 equilibria of the form $\{ (x,y) \ | \ x,y \in \{ R,S,P \} \}$, where $R=(1,0),S=(0,1),P=(0,0)$, which are all of (2,2) saddle type.
\end{prop5*}
\subsection{Nash equilibrium}
\begin{df}
$(x^*,y^*)$ is a Nash equilibrium of the system iff $(x^*)^T A y^* > x A y^*$ and $(y^*)^T B x^* > y^T B x^*$, for all $x, y$. 
\end{df}
In the bimatrix game in which there is only one Nash equilibrium, it can be found as the intersection of the nullclines (here of the equation (\ref{rownanie})). 
It is commonly known that system (\ref{rownanie}), is Hamiltonian for $A=-B^T$, that is when $\epsilon_X+ \epsilon_Y = 0$ (see \cite{[9]}, \cite{[16]}). Moreover it was proved that in this case all trajectories have on average the same payoff as the Nash equilibrium \cite{[18]}. On the other hand, on any given step huge deviations from the payoff, when playing Nash equilibrium, can be observed. Thus, a risk averse agent would prefer the Nash equilibrium to a chaotic orbit. \cite{[16]}
\\ \\
In the sixth section we will make use of the following proposition.
\begin{prop}\label{NashEq}
Nash equilibrium of the system  $(x^*,y^*)=(\frac{1}{3},\frac{1}{3},\frac{1}{3},\frac{1}{3})$, which is of $(2,2)$ saddle type for $\epsilon_X + \epsilon_Y \neq 0$ (or a center in the Hamiltonian case), with eigenvalues 
$$(\lambda_1,\lambda_2,\lambda_3,\lambda_4)= \Big(- \frac{1}{3} \sqrt{-3 + \epsilon_X \epsilon_Y - \sqrt{3} \sqrt{-(\epsilon_X + \epsilon_Y)^2}}, \frac{1}{3} \sqrt{-3 + \epsilon_X \epsilon_Y - \sqrt{3} \sqrt{-(\epsilon_X + \epsilon_Y)^2}}, 
$$
$$
- \frac{1}{3} \sqrt{-3 + \epsilon_X \epsilon_Y + \sqrt{3} \sqrt{-(\epsilon_X + \epsilon_Y)^2}}, \frac{1}{3} \sqrt{-3 + \epsilon_X \epsilon_Y + \sqrt{3} \sqrt{-(\epsilon_X + \epsilon_Y)^2}} \  \Big), 
$$
where $Re(\lambda_1),Re(\lambda_3) \leq 0$, $Re(\lambda_2),Re(\lambda_4) \geq 0$
\begin{proof}
One can see that 
$$ 
Ay= \begin{bmatrix}
-1 +\epsilon_X (- \frac{1}{3} +y_1) + y_1 +2y_2 \\
1- 2y_1 +\epsilon_X(- \frac{1}{3}+ y_2) -y_2 \\
y_1 +\epsilon_X( \frac{2}{3} -y_1 -y_2) -y_2 
\end{bmatrix} \  and \ \ 
(x^T B)^T= \begin{bmatrix}
1 +\epsilon_Y (- \frac{1}{3} +x_1) - x_1 -2x_2 \\
-1+ 2x_1 +\epsilon_Y (- \frac{1}{3}+ x_2) +x_2 \\
-x_1 +\epsilon_Y( \frac{2}{3} -x_1 -x_2) +x_2 
\end{bmatrix}  
$$
in particular 
\newline
$A[ \frac{1}{3},\frac{1}{3},\frac{1}{3}]^T= [ \ 0 \ \ 0 \ \ 0 \ ]^T$  and \ \ 
$[\frac{1}{3},\frac{1}{3},\frac{1}{3}] B= [ \ 0 \ \ 0 \ \ 0 \ ]$
which by the definition means that $(\frac{1}{3},\frac{1}{3},\frac{1}{3},\frac{1}{3})$ is a Nash equilibrium. Note that for $z_1:=-3 +\epsilon_X \epsilon_Y - \sqrt{3} \sqrt{-(\epsilon_X + \epsilon_Y)^2}$ and $z_2:=-3 +\epsilon_X \epsilon_Y + \sqrt{3} \sqrt{-(\epsilon_X + \epsilon_Y)^2}$, we have $Re(z_1)<0, Im(z_1) \leq 0$, so it follows that $Re(\sqrt{z_1}) \leq 0$. On the other hand $Re(z_2)<0, Im(z_2) \geq 0$, implies $Re(\sqrt{z_2}) \geq 0$. Hence $Re(\lambda_1),Re(\lambda_3) \leq 0$, $Re(\lambda_2),Re(\lambda_4) \geq 0$. Moreover we see that: 
$Re(\lambda_1)=0 \Leftrightarrow Re(z_1)<0 \Leftrightarrow Im(z_1)=0 \Leftrightarrow Re(\lambda_1)=Re(\lambda_3)=Re(\lambda_2)=Re(\lambda_4)=0 \Leftrightarrow Re(z_1)<0 \Leftrightarrow Im(z_1)=0 \Leftrightarrow \epsilon_X +\epsilon_Y =0$. In such a case we have that $B=-A^T$. Consider then a Lyapunov function 
$$
V(x_1,x_2,x_3,y_1,y_2,y_3):=\frac{1}{3} \sum{\log x_i} + \frac{1}{3} \sum{\log y_i}
$$
From Jensen's inequality $V(x_1,x_2,x_3,y_1,y_2,y_3) \leq 2 \log \frac{1}{3}$, with equality if and only if $x_1=x_2=x_3=y_1=y_2=y_3=\frac{1}{3}$. Moreover since $\sum{(Ay)_i}= \sum{(Bx)_i}=0$ and $B=-A^T$:
$$\dot{V}(x_1,x_2,x_3,y_1,y_2,y_3)= \frac{1}{3} \sum{(Ay)_i} - x^T Ay + \frac{1}{3} \sum{(Bx)_i} - y^T Bx= -x^T Ay + y^T A^T x= 
0$$ 
So either $(\frac{1}{3},\frac{1}{3},\frac{1}{3},\frac{1}{3})$ is a hyperbolic point of $(2,2)$ type or is a center. 
\end{proof}
\end{prop}

\subsection{Dynamics on the boundaries of codimensions 2 and 1}
If each of the players is restricted to only two strategies, then all orbits are periodic and spiral around one of the 6 equilibria or they converge to draw $(\{R,R\},\{P,P\},\{S,S\})$ with time tending to $\pm \infty$. The latter case occurs for the system restricted to the subspaces $\{x_i=0, \ y_i=0 \}$, $i=1,2,3$.
\begin{prop}\label{6addEq}
There are 6 equilibria $Z^a,Z^b,Z^c,Z^d,Z^e,Z^f$ of (1,2,1) type (1-dim. stable manifold, 2-dim. center manifold, 1-dim. unstable manifold), which are centers within the 2 dimensional invariant subspaces $H^a,H^b,H^c,H^d,H^e,H^f$ (where the system is integrable), and have one dimensional stable and unstable manifolds $W^s_a,W^u_a,W^s_b,W^u_b$, $W^s_c,W^u_c,W^s_d,W^u_d,W^s_e,W^u_e,W^s_f,W^u_f$: 
\newline
$$Z^a= \Big(0, \frac{2}{3-\epsilon_Y},\frac{1+\epsilon_X}{3+\epsilon_X},\frac{2}{3+\epsilon_X}\Big), H^a= \{x_1=0,y_3=0 \}, W^s_a \subset \{ x_1=0 \}, W^u_a \subset \{ y_3=0 \}$$

$$Z^b= \Big( \ 0, \ \frac{1+\epsilon_Y}{3+\epsilon_Y}, \ \frac{1-\epsilon_X}{3-\epsilon_X}, \ 0 \ \Big), \  H^b =\{x_1=0,y_2=0\}, \  W^s_b \subset \{ y_2=0 \}, \ W^u_b \subset \{ x_1=0 \}$$  

$$Z^c= \Big( \ \frac{1-\epsilon_Y}{3-\epsilon_Y}, \ 0 , \ 0, \ \frac{1+\epsilon_X}{3+\epsilon_X} \Big), \  H^c= \{x_2=0,y_1=0 \}, \ W^s_c \subset \{ x_2=0 \}, \ W^u_c \subset \{ y_1=0 \}$$

$$Z^d= \Big(\frac{2}{3+\epsilon_Y}, 0,\frac{2}{3-\epsilon_X},\frac{1-\epsilon_X}{3-\epsilon_X}\Big), H^d = \{ x_2=0,y_3=0 \}, W^s_d \subset \{ y_3=0 \}, W^u_d \subset \{ x_2=0 \}$$

$$Z^e= \Big(\frac{1+\epsilon_Y}{3+\epsilon_Y}, \frac{2}{3+\epsilon_Y},0,\frac{2}{3-\epsilon_X}\Big), H^e = \{x_3=0,y_1=0\}, W^s_e \subset \{ y_1=0 \}, W^u_e \subset \{ x_3=0 \}$$

$$Z^f= \Big( \frac{2}{3-\epsilon_Y}, \frac{1-\epsilon_Y}{3-\epsilon_Y},\frac{2}{3+\epsilon_X},0 \Big), H^f = \{x_3=0,y_2=0 \}, W^s_f \subset \{ x_3=0 \}, W^u_f \subset \{ y_2=0 \}$$

\begin{proof}
For $Z^a$ one can compute the eigenvalues of linearized system to be 
$$
\{ \frac{3+{\epsilon_X}^2}{3+\epsilon_X}, - 2i \sqrt{\frac{1+\epsilon_X}{3+\epsilon_X}}\sqrt{\frac{1-\epsilon_Y}{3-\epsilon_Y}}, 2i \sqrt{\frac{1+\epsilon_X}{3+\epsilon_X}}\sqrt{\frac{1-\epsilon_Y}{3-\epsilon_Y}}, -\frac{3+{\epsilon_Y}^2}{3-\epsilon_Y} \}
$$ 
Eigenvectors corresponding to purely imaginary eigenvalues lie within subspace $x_1=0$, $y_1+y_2=1$. Restricting our system to this subspace (that is substituting $y_2 = 1-y_1$ and $x_1=0$ in the equation), we find the Lyapunov function for $Z^a$ within subspace $x_1=0$, $y_1+y_2=1$: 
$$
V_a(x_2,y_1) := (1-\epsilon_Y) \log (1-x_2) + 2 \log x_2 + (1+ \epsilon_X) \log y_1 + 2 \log (1-y_1)
$$
One can see that $V_a(x_2,y_1) \leq 0$ for $x_2,y_1 \in (0,1)$ (moreover $V_a(x_2,y_1)=0 \Leftrightarrow x_2=\frac{2}{3-\epsilon_Y}, y_1=\frac{1+\epsilon_X}{3+\epsilon_X}$) and 
\begin{multline*}
\dot{V_a}(x_2,y_1)= (1-\epsilon_Y) \cdot (\frac{-x_2 '}{1-x_2}) + 2 \frac{x_2 '}{x_2} + (1+ \epsilon_X) (\frac{y_1 '}{y_1}) + 2 \frac{-y_1 '}{1-y_1} = (1-\epsilon_Y) \cdot x_2 (-1+\epsilon_X (-1+y_1)+3 y_1) + 
\\ 
2 (-1+ x_2) (-1+\epsilon_X (-1+y_1)+3 y_1) + (1+ \epsilon_X) (2 + (-3 + \epsilon_Y) x_2) (-1 + y_1) + 2 (2 + (-3 + \epsilon_Y) x_2) y_1=0
\end{multline*}
which proves that $Z^a$ is a center within subspace $x_1=0$, $y_1+y_2=1$.
\newline
To prove that, for $\epsilon_X+\epsilon_Y=0$, stable manifold of $Z^a$ is contained in the set $\{ x_1=0\}$, let us note that the system is Hamiltonian in $int (\Delta \times \Delta)$ so if some part of this stable manifold was included in $int (\Delta \times \Delta)$ then on the whole of unstable manifold the function $V$, defined in the previous proposition, would be constant which is impossible since $V \rightarrow -\infty$ when $x_1 \rightarrow 0$.
\newline
In general, note that the eigenvector corresponding to the eigenvalue $\frac{3+{\epsilon_Y}^2}{-3+\epsilon_Y}$ is equal to 
\begin{multline*}
v=(0,v_2,v_3,1)=
\\
\Big( 0, \frac{(-1 + \epsilon_X) (3 + \epsilon_X)^2 (-1 + \epsilon_Y) (3 + \epsilon_Y^2)}{(3- \epsilon_Y) (-3 (-5 + \epsilon_X) (3 + 2 \epsilon_X) + 8 (3 - \epsilon_X) (1 + \epsilon_X) \epsilon_Y + 2 (-12 + (-5 + \epsilon_X) \epsilon_X) \epsilon_Y^2 + (3 + \epsilon_X) \epsilon_Y^4)}, 
\\
-\frac{(1 + \epsilon_X) (3 (3 + \epsilon_Y^2)^2 + \epsilon_X (33 + \epsilon_Y (-32 + 14 \epsilon_Y + \epsilon_Y^3)))}{6 (-5 + \epsilon_X) (3 + 2 \epsilon_X) - 16 (-3 + \epsilon_X) (1 + \epsilon_X) \epsilon_Y + 4 (-12 + (-5 + \epsilon_X) \epsilon_X) \epsilon_Y^2 - 2 (3 + \epsilon_X) \epsilon_Y^4}, 1 \Big)
\end{multline*}
Let us write 
$$x_1=0+h_1, x_2= \frac{2}{3-\epsilon_Y} +h_2, y_1= \frac{1+\epsilon_X}{3+\epsilon_X}+s_1, y_2=\frac{2}{3+\epsilon_X}+s_2
$$ 
It is easy to see that for all $\epsilon_X, \epsilon_Y \in (-1,1)$ we have $v_2,v_3 \in (0,1)$. So for a point $(x_1,x_2,y_1,y_2)$ which belongs to the stable manifold of the point $Z^a$ and is very close to $Z^a$, we may assume that 
$$
sign(h_2)=sign(v_2)>0, \ sign(s_1)=sign(v_3)>0, \ sign(s_2)=sign(1)>0
$$ 
(the second case when $sign(h_2)=-sign(v_2)<0$, $sign(s_1)=-sign(v_3)<0$, $sign(s_2)=-sign(1)<0$ can be treated in the same way). We want to prove that projection of the stable manifold on the $X$ -coordinate lies in the set $x_1=0$, so it suffices to prove that with the above assumptions on $h_2,s_1,s_2$ and with $h_1=0$ we have $h_2 ' <0$ (or  $h_2 ' >0$ if we change the signs of $h_2,s_1,s_2$). Obviously $h_1 ' = 0$ on the subset $h_1=0$ and 
$$
x_2 ' = h_2 ' = \big( 1 + \epsilon_Y (-1+ h_2)- 3 h_2 \big)  \big( 2 + (3 - \epsilon_Y) h_2 \big) \frac{(-3 + \epsilon_X) s_1 + 2 \epsilon_X s_2}{(-3 + \epsilon_Y)^2}
$$ 
Since $h_2>0$ and is small enough then 
$$
(1 + \epsilon_Y (-1 + h_2)-3 h_2) (2 + (3 - \epsilon_Y) h_2)>0
$$
If $\epsilon_X \leq0$ then obviously 
$$
(-3 + \epsilon_X) s_1 + 2 \epsilon_X s_2<0
$$ since $s_1,s_2>0$, otherwise we have 
$$
(-3 + \epsilon_X) s_1 + 2 \epsilon_X s_2<0 \Leftrightarrow  \frac{2 \epsilon_X}{3-\epsilon_X} < \frac{s_1}{s_2}
$$
Since for any fixed $\epsilon_X, \epsilon_Y \in (-1,1)$ the following inequality is satisfied $\frac{2 \epsilon_X}{3-\epsilon_X} < v_3$, we can choose $s_1,s_2$ so small such that $\frac{s_1}{s_2}$ is close enough to $v_3$ and then $\frac{2 \epsilon_X}{3-\epsilon_X} < \frac{s_1}{s_2}$, which proves that projection of stable manifold on the $X$ -coordinate lies in the set $x_1=0$
\end{proof}
\end{prop}

\begin{figure}[!htb]
\minipage{0.23\textwidth}
  \includegraphics[width=\linewidth]{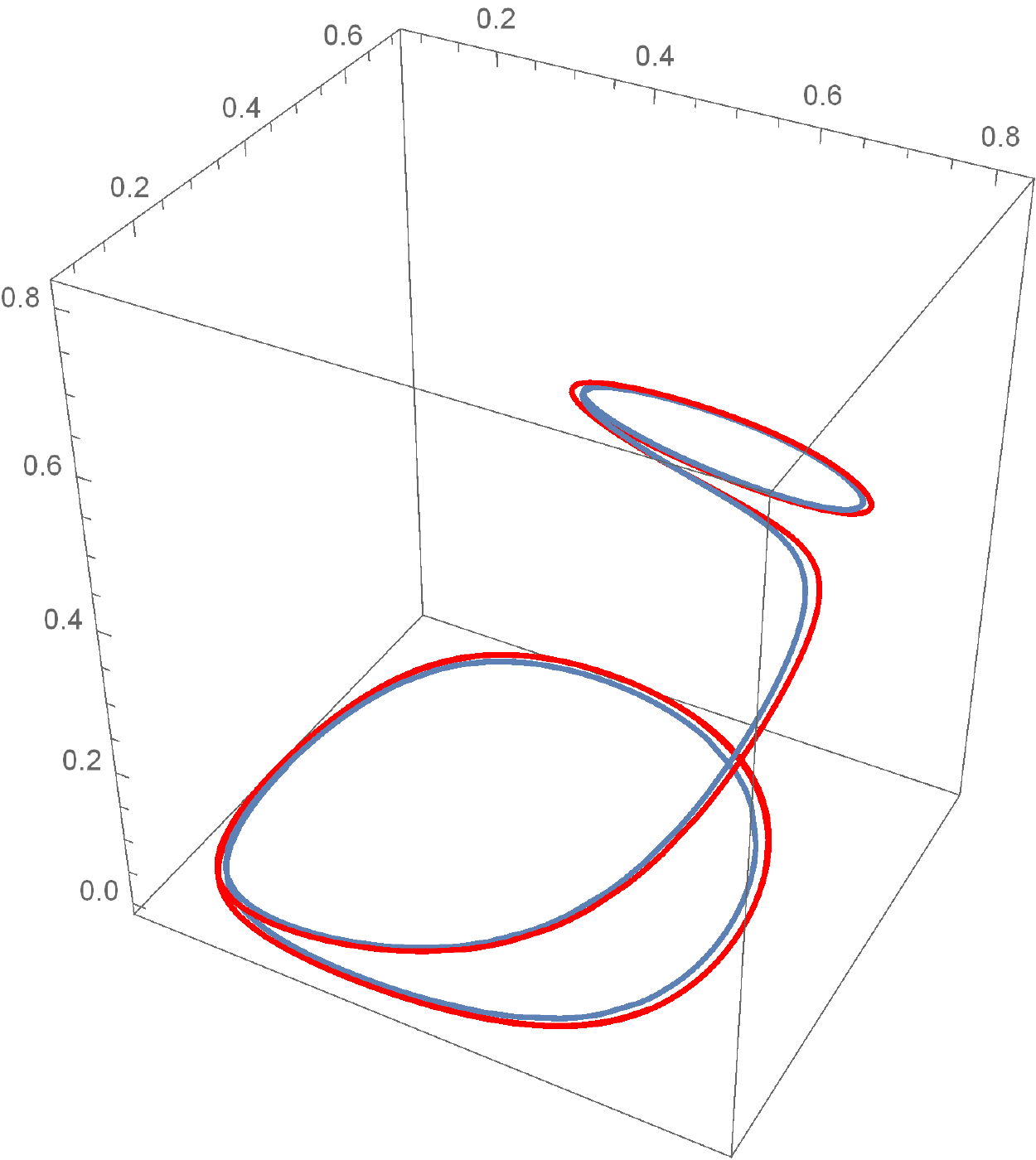}
	{\tiny (a) $\epsilon_X= -0.09, \ \epsilon_Y= -0.79$}
\endminipage\hfill
\minipage{0.23\textwidth}
  \includegraphics[width=\linewidth]{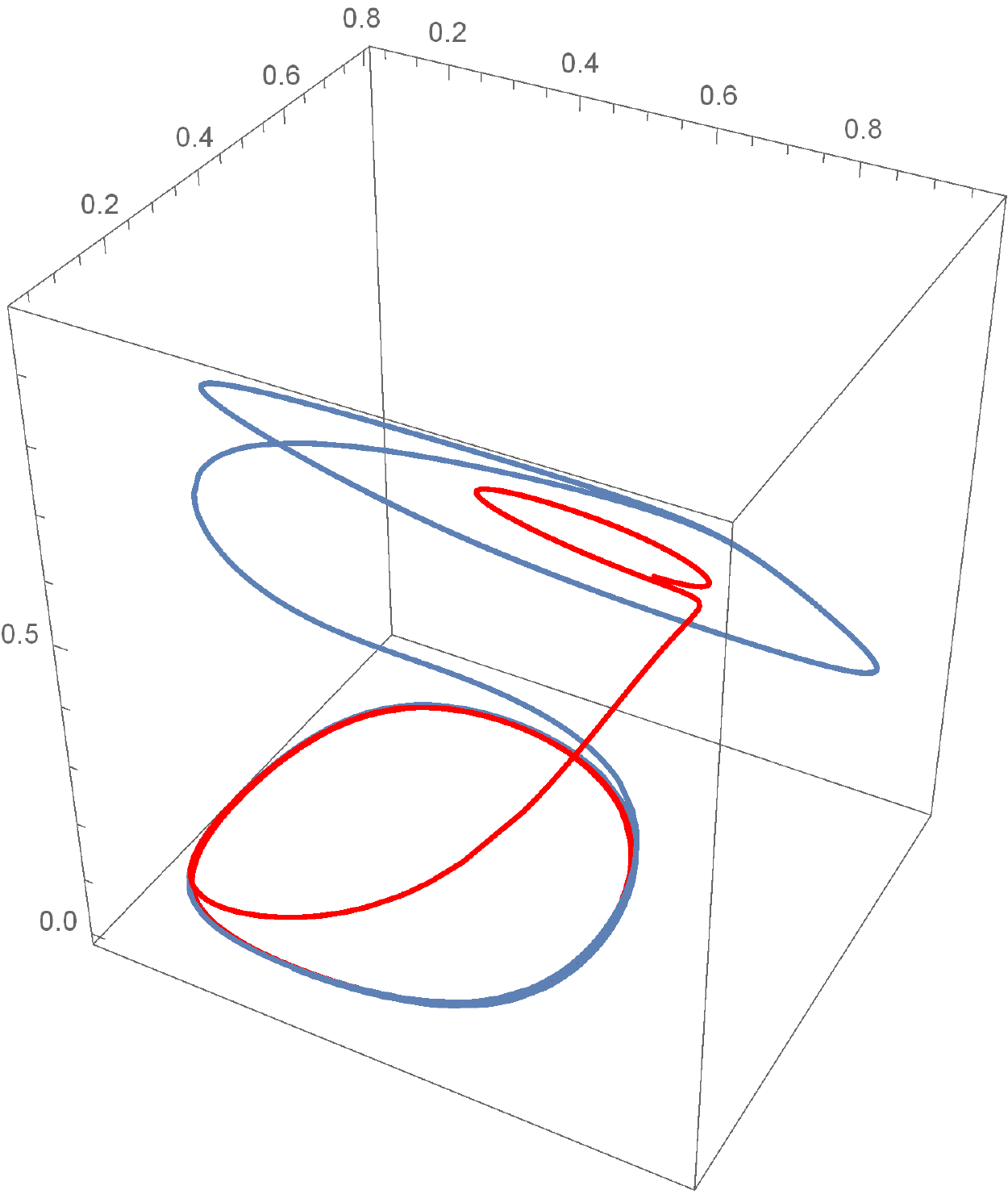}
	{\tiny (b) $\epsilon_X= 0.81, \ \epsilon_Y= 0.11$}
\endminipage\hfill
\minipage{0.23\textwidth}%
  \includegraphics[width=\linewidth]{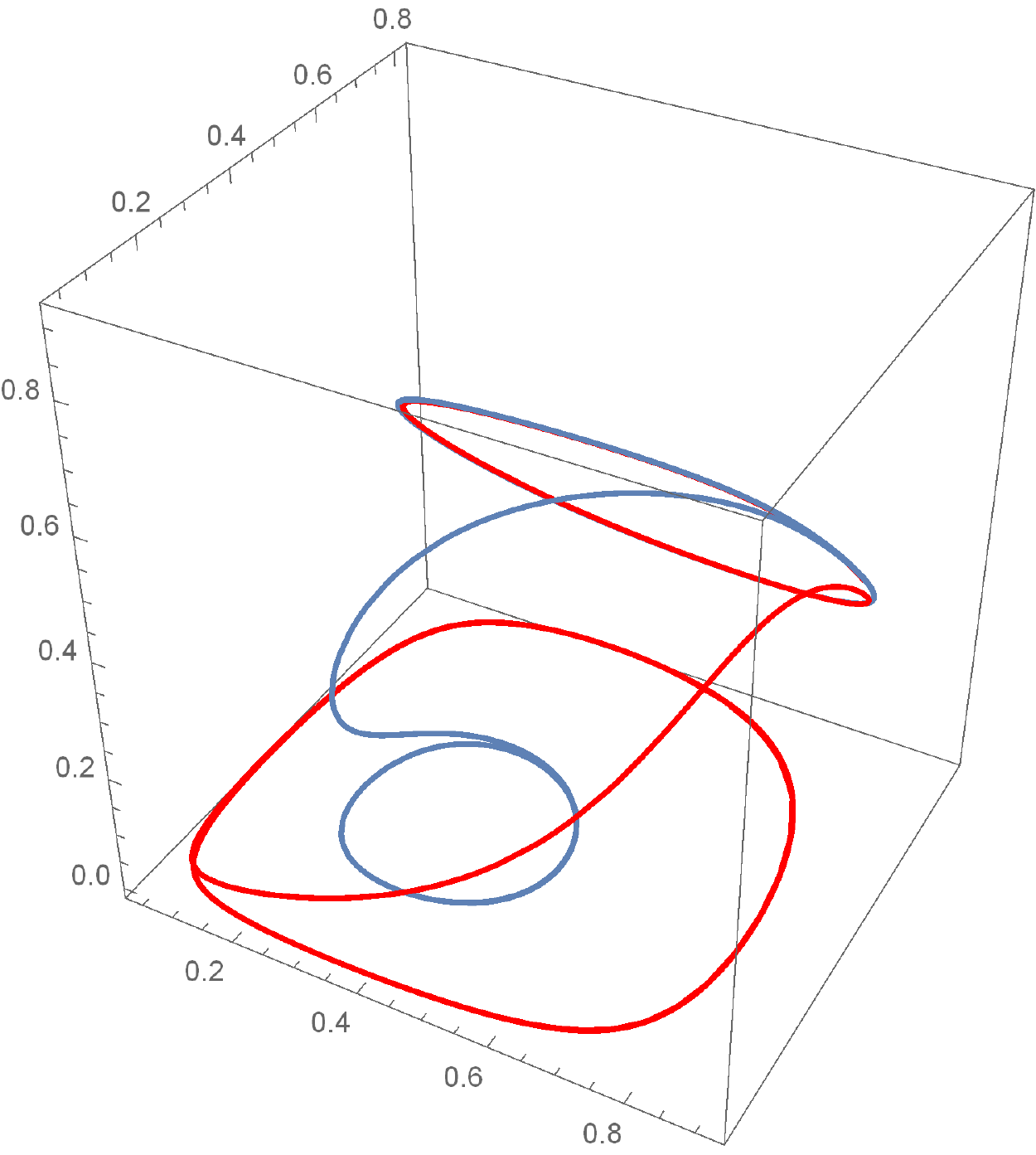}
	{\tiny (c) $\epsilon_X= 0.41, \ \epsilon_Y= 0.81$}
\endminipage\hfill
\minipage{0.23\textwidth}%
  \includegraphics[width=\linewidth]{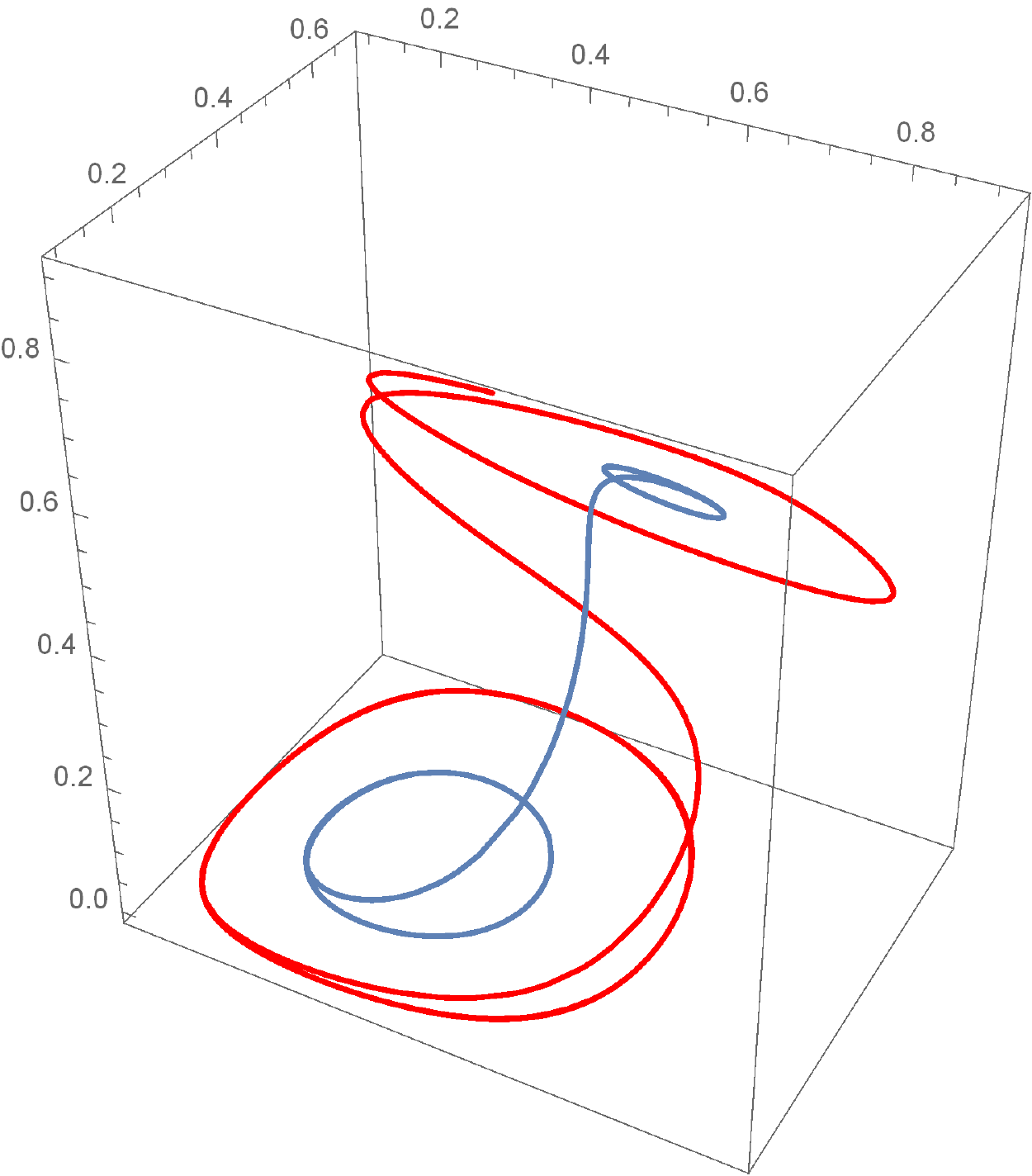}
	{\tiny (d) $\epsilon_X= -0.79, \ \epsilon_Y= -0.29$}
\endminipage
\caption{Generic types of behaviour of the flow on the boundary $\{ x_1=0 \}$, red curve is a a backward orbit of starting point close to the stable manifold of $Z^a$, blue curve is a forward orbit of a starting point close to the unstable manifold of $Z^b$. The curves stay either (a) close to each other for all times, (b) close near the subspace $\{ y_2=0 \}$, (c) close near $\{ y_3=0 \}$, (d) or are faraway from each other.}
\end{figure}
On the boundary of codimension 1, the trajectories of the game travel from one invariant hyperplane of dimension 2 to another one.
\begin{prop}
For all $\epsilon_X,\epsilon_Y \in (-1,1)$, the following inclusions hold:
$$\alpha(W^{s}_{Z^a}) \subset \{ x_1=0, \ y_2=0 \}, \  \alpha(W^{s}_{Z^b}) \subset \{ x_3=0, \ y_2=0 \}, \  \alpha(W^{s}_{Z^c}) \subset \{ x_2=0, \ y_3=0 \}$$
$$\alpha(W^{s}_{Z^d}) \subset \{ x_1=0, \ y_3=0 \}, \ \alpha(W^{s}_{Z^e}) \subset \{ x_2=0, \ y_1=0 \},  \  \alpha(W^{s}_{Z^f}) \subset \{ x_3=0, \ y_1=0 \}$$
$$\omega(W^{u}_{Z^a}) \subset \{ x_2=0, \ y_3=0 \}, \ \omega(W^{u}_{Z^b}) \subset \{ x_1=0, \ y_3=0 \}, \ \omega(W^{u}_{Z^c}) \subset \{ x_3=0, \ y_1=0 \}$$
$$\omega(W^{u}_{Z^d}) \subset \{ x_2=0, \ y_1=0 \}, \ \omega(W^{u}_{Z^e}) \subset \{ x_3=0, \ y_2=0 \}, \ \omega(W^{u}_{Z^f}) \subset \{ x_1=0, \ y_2=0 \}$$
\begin{proof}
For $W^{s}_{Z^a}$ and $W^{u}_{Z^b}$, note that with $x_1=0$ in the primary equation (\ref{rownanie}) we have:
\newline
$\dot{y_2}>0$ if $x_2 \in (0, \frac{1}{2}]$ and $\dot{y_3}<0$ if $x_2 \in [\frac{1}{2} ,1)$. 
Moreover $\dot{x_2}>0$ holds when either $(\epsilon_X<0$ and $0<y_1<\frac{1+\epsilon_X}{3+\epsilon_X})$ or $(\epsilon_X>0$ and $0<y_1<\frac{1-\epsilon_X}{3-\epsilon_X})$. 
Since $\frac{2}{3-\epsilon_Y}>\frac{1}{2}$ and $\frac{1+\epsilon_Y}{3+\epsilon_Y}<\frac{1}{2}$, this along with the above, means that the forward orbit of a point from $W^{u}_{Z^b}$ will accumulate on the hyperplane $\{ x_1=0, \ y_3=0 \}$ and the backward orbit of a point from $W^{s}_{Z^a}$ will accumulate on the hyperplane $\{ x_1=0, \ y_2=0 \}$, where the motion is periodic.
\end{proof}
\end{prop}

We believe that the following conjecture is worth further investigation.
\begin{conj}
As one can check, there exist parameter values (e.g. $\epsilon_X= -0.09$, $\epsilon_Y= -0.79$) for which stable manifold of $Z^a$ and unstable of $Z^b$ seem to be very close. 
\newline
This observation motivates the following question:
\newline
The set of pairs $(\epsilon_X,\epsilon_Y)$ for which there is a heteroclinic orbit from $Z^b$ to $Z^a$ is isolated in the parameter space. The same question remains for the existence of other heteroclinic orbits, that is from:
\newline
$Z^a$ to $Z^d$, $Z^d$ to $Z^c$, $Z^c$ to $Z^e$, $Z^e$ to $Z^f$, $Z^f$ to $Z^b$
\end{conj}

\begin{rmq}
The existence of the heteroclinic orbit from $Z^b$ to $Z^a$ is equivalent to the existence of the zero of the function:
$$
(s,t,\epsilon_X,\epsilon_Y) \longmapsto g_b (s,\epsilon_X,\epsilon_Y) - g_a(t,\epsilon_X,\epsilon_Y)
$$
where $g_a, g_b$ are parameterizations (by s and t, respectively) of the manifolds $W^{s}_{Z^a}$ and $W^{u}_{Z^b}$. This is a codimension 2 problem, hence in a family of two parameter vector fields one should expect the existence of the solution for the isolated set of parameters. (see \cite{[22]}, \cite{[27]}, \cite{[28]}, \cite{[29]})
\end{rmq}

\section{Heteroclinic network and its transition maps}
Here we will show how heteroclinic network naturally appears in our system. We describe all of the maps associated with this network and process of deriving it. We will indicate the differences between our and Aguiar\&Castro model. Although that one was incorrect, we have to take their observation into account:
\begin{rmq}
One can check that the system is $C^1$ linearizable for all $\epsilon_X,\epsilon_Y \in (-1,1)$ except the case $\epsilon_X=\epsilon_Y$. The whole forthcoming analysis is conducted for $\epsilon_X \neq \epsilon_Y$.
\end{rmq}
\begin{prop}\label{9Equi}
The equation (\ref{rownanieZredukowane}), possesses 9 equilibria of the form $\{ (x,y) \ | \ x,y \in \{ R,S,P \} \}$, where $R=(1,0),S=(0,1),P=(0,0)$. They are all of (2,2) saddle type with eigenvalues (and corresponding eigenvectors) and stable/unstable manifolds $W^s/W^u$ 
$($since the system is invariant under the symmetry $\sigma$ it suffices to provide the description for points $(P,P)$, $(P,S)$, $(S,P)$$\ )$:
\begin{align*}
for \  (P,P), \  eigenbasis  \ (e_1,e_2,e_3,e_4): \  (\lambda_1, \lambda_2, \lambda_3, \lambda_4)= (-1-\epsilon_X, 1-\epsilon_X, -1-\epsilon_Y, 1-\epsilon_Y)
\\ 
W^s_{(P,P)}= \{ (x_1,0,y_1,0) \ | \ x_1 \neq 1, y_1 \neq 1  \} \cap \Delta, \ 
W^u_{(P,P)}= \{ (0,x_2,0,y_2) \ | x_2 \neq 1, y_2 \neq 1 \} \cap \Delta 
\end{align*}
\begin{align*}
for \ (S,P),  \ eigenbasis \ (e_1-e_2,e_2,e_3,e_4)):  (\lambda_1, \lambda_2, \lambda_3, \lambda_4)= (-2, -1+\epsilon_X, 2, 1+\epsilon_Y) 
\\
W^s_{(S,P)}= \{ (x_1,x_2,0,0) \} \cap \Delta \backslash \{ x_2=0 \}, \ W^u_{(S,P)}= \{ (0,1,y_1,y_2) \} \cap \Delta \backslash \{ y_1+y_2=1 \}
\end{align*}
\begin{align*}
for \ (P,S), \ eigenbasis \ (e_1,e_2,e_3-e_4,e_4):  (\lambda_1, \lambda_2, \lambda_3, \lambda_4)= (2, 1+\epsilon_X, 2, -2, -1+\epsilon_Y)
\\ 
W^s_{(P,S)} = \{ (0,0,y_1,y_2) \} \cap \Delta \backslash \{ y_2=0 \},  \ W^u_{(P,S)}= \{ (x_1,x_2,0,1) \} \cap \Delta \backslash \{ x_1+x_2=0 \}
\end{align*}
\begin{proof}
We find that the linearized system near $(P,S)$ is given by:
$$
\begin{bmatrix}
\dot{x_1} \\
\dot{x_2} \\
\dot{y_1} \\
\dot{y_2}
\end{bmatrix} = 
\begin{bmatrix}
2 && 0 						&& 0 						&& 0 \\
0	&& 1+\epsilon_X && 0 						&& 0 \\
0	&& 0 						&& -2 					&& 0 \\
0 && 0 						&& 1+\epsilon_Y && -1+\epsilon_Y
\end{bmatrix} \cdot
\begin{bmatrix}
x_1 \\
x_2 \\
y_1 \\
y_2 -1
\end{bmatrix} 
$$
Hence, the searched eigenvalues are $\{ 2, 1+\epsilon_X, -2, -1+\epsilon_Y \}$ with the corresponding eigenvectors 
\newline
$\{ (1,0,0,0), (0,1,0,0), (0,0,1,-1), (0,0,0,1) \}$. 
\newline
$W^s_{(P,S)}$ is tangent at $(P,S)$ to the 2-dimensional subspace $(P,S)+$ lin$\{ (0,0,1,-1), (0,0,0,1) \}$. Note that the subset $\{ x_1=0, x_2=0 \}$ is invariant, and moreover, when $y_1 \neq 0$, $y_2 \neq 0$, we have $\dot{y_1}<0$, $\dot{y_2}>0$. Furthermore the subset $\{ x_1=0, x_2=0, y_1=0 \}$ is invariant and $\dot{y_2}>0$ there. So the whole set $\{ (0,0,y_1,y_2) \ | \ y_2 \neq 0 \}$ is attracted to points $(P,S)= (0,0,0,1)$. From uniqueness of the existing stable manifold we have that $W^s_{(P,S)}= \{ (0,0,y_1,y_2) \} \cap \Delta \backslash \{ y_2=0 \}$. Analogically we prove $W^u_{(P,S)}= \{ (x_1,x_2,0,1) \} \cap \Delta \backslash \{ x_1+x_2=0 \}$
\end{proof}
\end{prop}
\begin{rmq}
There are no other equilibria for (\ref{rownanieZredukowane}) than described so far.
\end{rmq}
Let us denote $[T_1,T_2] := \{ (1- \lambda) T_1+ \lambda T_2 \ | \ \lambda \in (0,1) \}$. 
Due to the above proposition, there exist a heteroclinic network $\Sigma$ consisting of these 9 equilibria $\{ (x,y) \ | \ x,y \in \{ R,S,P \} \}$ and heteroclinic orbits written explicitly:
\newline 
$[(P,P),(P,S)]$, $[(P,P),(S,P)]$, $[(P,S),(S,S)]$, $[(P,S),(R,S)]$, $[(S,S),(S,R)]$, $[(S,S),(R,S)]$, 
\newline
$[(R,S),(R,R)]$, $[(R,S),(R,P)]$, $[(R,P),(P,P)]$, $[(R,P),(S,P)]$, $[(S,P),(S,S)]$, $[(S,P),(S,R)]$,
\newline
$((S,R),(R,R)]$, $[(S,R),(P,R)]$, $[(R,R),(R,P)]$, $[(R,R),(P,R)]$, $[(P,R),(P,P)]$, $[(P,R),(P,S)]$
\newline 
\begin{rmq}
Heteroclinic orbits in this network are non-transversal intersections of stable and unstable manifolds of the appropriate equilibria.
\end{rmq}
Because of the above remark and hyperbolicity of the equilibria from the network, one does not necessarily expect chaotic switching infinitely close to the heteroclinic network as indeed we show.
In fact we will discuss this matter more closely in sections 5 and 7.


 \tikzset{vertex/.style={shape=circle, 
                        minimum size=1pt,
                        fill=gray,
                        inner sep = 0pt}}  
\begin{figure}
\centering
\begin{tikzpicture}[xscale=1,>=latex]
       \coordinate  (PP) at (0,0);   \coordinate (RS) at (4,4);
       \coordinate  (RP) at (4,0);  \coordinate  (PS) at (0,4); 
       \coordinate  (SS) at (2,7.464);   \coordinate (RR) at (7.464,2); 
			 \coordinate  (SP) at (2,3.464); \coordinate (PR) at (3.464,2); 
			 \coordinate  (SR) at (5.3836,5.3836);   
\draw[yellow,decoration={markings, mark=at position 0.5 with {\arrow{>}}}, postaction={decorate}] (PP) -- (PS);
\draw[red,decoration={markings, mark=at position 0.5 with {\arrow{>}}}, postaction={decorate}] (PS) -- (RS);
\draw[red,decoration={markings, mark=at position 0.5 with {\arrow{>}}}, postaction={decorate}] (RS) -- (RP);
\draw[yellow,decoration={markings, mark=at position 0.5 with {\arrow{>}}}, postaction={decorate}] (RP) -- (PP);
\draw[yellow,decoration={markings, mark=at position 0.5 with {\arrow{>}}}, postaction={decorate}] (PS) -- (SS);
\draw[green,decoration={markings, mark=at position 0.5 with {\arrow{>}}}, postaction={decorate}] (SS) -- (RS);
\draw[green,decoration={markings, mark=at position 0.5 with {\arrow{>}}}, postaction={decorate}] (RS) -- (RR);
\draw[yellow,decoration={markings, mark=at position 0.5 with {\arrow{>}}}, postaction={decorate}] (RR) -- (RP);
\draw[green,decoration={markings, mark=at position 0.5 with {\arrow{>}}}, postaction={decorate}] (PP) -- (SP);
\draw[green,decoration={markings, mark=at position 0.5 with {\arrow{>}}}, postaction={decorate}] (SP) -- (SS);
\draw[red,decoration={markings, mark=at position 0.5 with {\arrow{>}}}, postaction={decorate}] (SP) -- (SR);
\draw[yellow,decoration={markings, mark=at position 0.5 with {\arrow{>}}}, postaction={decorate}] (SS) -- (SR);
\draw[yellow,decoration={markings, mark=at position 0.5 with {\arrow{>}}}, postaction={decorate}] (SR) -- (RR);
\draw[red,decoration={markings, mark=at position 0.5 with {\arrow{>}}}, postaction={decorate}] (SR) -- (PR);
\draw[red,decoration={markings, mark=at position 0.5 with {\arrow{>}}}, postaction={decorate}] (RP) -- (SP);
\draw[green,decoration={markings, mark=at position 0.5 with {\arrow{>}}}, postaction={decorate}] (RR) -- (PR);
\draw[green,decoration={markings, mark=at position 0.5 with {\arrow{>}}}, postaction={decorate}] (PR) -- (PP);
\draw[red,decoration={markings, mark=at position 0.5 with {\arrow{>}}}, postaction={decorate}] (PR) -- (PS);
\node at (-0.4,-0.2) {\scriptsize $(P,P)$};
\node at (4.2,-0.2) {\scriptsize $(R,P)$};
\node at (-0.4,4.2) {\scriptsize $(P,S)$};
\node at (3.6,3.8) {\scriptsize $(R,S)$};
\node at (1.6,7.6) {\scriptsize $(S,S)$};
\node at (5.8,5.5) {\scriptsize $(S,R)$};
\node at (8,2) {\scriptsize $(R,R)$};
\node at (3.6,1.75) {\scriptsize $(P,R)$};
\node at (1.55,3.5) {\scriptsize $(S,P)$};
\end{tikzpicture}
\caption{RSP network}
\end{figure}
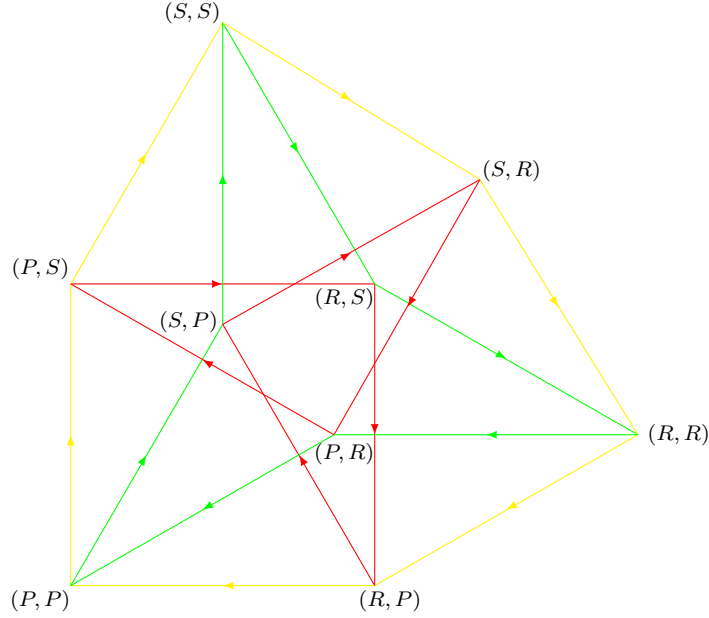

\begin{prop}
The following matrices $MacVU$ $($for $V,U \in \{ R,P,S \} )$ are the transformation matrices from the eigenbasis at the point $(V,U)$ to the canonical basis $(e_1,e_2,e_3,e_4)$ of $\mathbb R^4$:
\newline
$$
MacRR = [e_1,e_2-e_1,e_3,e_4-e_3] = \left( \begin{array}{llll}
1 & -1 & 0 & 0 \\
0 & 1 & 0 & 0 \\
0 & 0 & 1 & -1 \\
0 & 0 & 0 & 1 
\end{array} \right)
$$
\begin{flushleft}
$MacRS = [e_1,e_2-e_1,e_3-e_4,e_4], \ 
MacSP = [e_1-e_2,e_2,e_3-e_4,e_4], \ 
MacPP = [e_1,e_2,e_3,e_4]$
\end{flushleft}
\begin{flushleft}
$MacRP = [e_1,e_2-e_1,e_3,e_4], \ 
MacPR \ = \ [e_1,e_2,e_3,e_4-e_3], \
MacPS = [e_1,e_2,e_3-e_4,e_4]$
\end{flushleft}
\begin{flushleft}
$MacSR \ = \ [e_1-e_2,e_2,e_3,e_4-e_3], \ 
MacSS \ = \ [e_1-e_2,e_2,e_3-e_4,e_4]$
\end{flushleft}
\begin{proof}
$MacVU$ consists of the eigenvectors of the linearized system at the point $(V,U)$.
\end{proof}
\end{prop}

\begin{prop} 
At each of the points $(V,U)$ $($with $V,U \in \{ R,P,S \})$ from the heteroclinic network, the local coordinates $(e'_1,e'_2,e'_3,e'_4)$ of affine space $\mathbb R^4$ in which this point is an origin, and the linearized system at this point is represented by the diagonal matrix (consisting of eigenvalues), are given by:
$$
(e_1',e_2',e_3',e_4') \ = \ MacVU^{-1} \big( (e_1, e_2, e_3, e_4) - (V,U) \big)
$$ 
\begin{proof}
Note that the transformation matrix from the eigenbasis (that is $(e_1, e_2, e_3 - e_4, e_4)$) at the point $(P,S)$ to the canonical basis of $\mathbb R^4$ is given by transformation matrix: 
$$
MacPS = [ e_1, e_2, e_3 - e_4, e_4 ] =  
\begin{bmatrix}
1 && 0 && 0 && 0 \\
0 && 1 && 0 && 0 \\
0 && 0 && 1 && 0 \\
0 && 0 && -1 && 1 
\end{bmatrix}
$$
Moreover one can check that MacVU$(V,U)=(V,U)$, hence the coordinates $(e_1',e_2',e_3',e_4')$ satisfying statements above should be given by condition: $MacPS ((e_1',e_2',e_3',e_4') + (P,S)) = (e_1,e_2,e_3,e_4)$. This gives the thesis.
\end{proof}
\end{prop}
Near each of the equilibrium from the heteroclinic network, we set four cross-sections to catch the orbits passing by close to this equilibrium. By introducing these cross-sections we can derive local maps near equilibriums and global maps from the neighbourhood of one equilibrium to another one along heteroclinic orbit, and in consequence discretize system in the neighbourhood of the heteroclinic network. Hence we can also check what kind of switching is possible and whether some particular trajectories infinitely close to the network are attainable.

\begin{lm}\label{sections}
The in/out sections near each of the points from the heteroclinic network, associated with the heteroclinic orbits between two points might be chosen as follows ($h>0$ is arbitrarily small):
\\ \\
In $(P,P)$ local coordinates:
$$
\Sigma^{out,(P,P)}_{(P,P) \rightarrow (P,S)} = \{ (z_1, z_2, z_3, h) \ | \ |z_1|,|z_2|,|z_3| < h \}, \ 
\Sigma^{in,(P,P)}_{(P,P) \rightarrow (P,S)} = \{ (z_1, z_2, z_3, 1-h) \ | \ |z_1|,|z_2|,|z_3| < h \},
$$
$$
\Sigma^{out,(P,P)}_{(P,P) \rightarrow (S,P)} = \{ (z_1, h, z_3, z_4) \ | \ |z_1|,|z_3|,|z_4| < h \}, \ 
\Sigma^{in,(P,P)}_{(P,P) \rightarrow (S,P)} = \{ (z_1, 1 -h, z_3, z_4) \ | \ |z_1|,|z_3|,|z_4| < h \}$$
\newline
In $(P,S)$ local coordinates:
$$
\Sigma^{out,(P,S)}_{(P,S) \rightarrow (R,S)} = \{ ( h, z_2, z_3, z_4) \ | \ |z_2|,|z_3|,|z_4| < h \}, \ 
\Sigma^{in,(P,S)}_{(P,S) \rightarrow (R,S)} = \{ (1-h, z_2, z_3, z_4) \ | \ |z_2|,|z_3|,|z_4| < h \},
$$
$$
\Sigma^{out,(P,S)}_{(P,S) \rightarrow (S,S)} = \{ (z_1, h, z_3, z_4) \ | \ |z_1|,|z_3|,|z_4| < h \}, \ 
\Sigma^{in,(P,S)}_{(P,S) \rightarrow (S,S)} = \{ (z_1, 1 - h, z_3, z_4) \ | \ |z_1|,|z_3|,|z_4| < h \}
$$
\newline
In $(S,P)$ local coordinates:
$$
\Sigma^{out,(S,P)}_{(S,P) \rightarrow (S,S)} = \{ ( z_1, z_2, z_3, h) \ | \ |z_1|,|z_2|,|z_3| < h \}, \ 
\Sigma^{in,(S,P)}_{(S,P) \rightarrow (S,S)} = \{ (z_1, z_2, z_3, 1-h) \ | \ |z_1|,|z_2|,|z_3| < h \},
$$
$$
\Sigma^{out,(S,P)}_{(S,P) \rightarrow (S,R)} = \{ (z_1, z_2, h, z_4) \ | \ |z_1|,|z_2|,|z_4| < h \}, \ 
\Sigma^{in,(S,P)}_{(S,P) \rightarrow (S,R)} = \{ (z_1, z_2, 1 -h, z_4) \ | \ |z_1|,|z_2|,|z_4| < h \}$$
\begin{proof}
Recall that in $(R,R)$ coordinates we have: $(R,R)=(0,0,0,0)$ and $(R,P)=MacRR^{-1} \big( (1,0,0,0)-(1,0,1,0) \big)=(0,0,-1,0)$ hence the heteroclinic orbit from $(R,R)$ to $(R,P)$ is parameterized by $(0,0,-\lambda,0)$ for $\lambda \in (0,1)$. So in $(R,R)$ coordinates, we define the "out" section $\Sigma^{out,(R,R)}_{(R,R) \rightarrow (R,P)}$ from point $(R,R)$ towards point $(R,P)$, to be the set $\{ (z_1,z_2,-h,z_4) \ | \ |z_1|, |z_2|, |z_4| <h \}$ for $h>0$ choosen arbitrarily small. 
Analogically, we define the "in" section $\Sigma^{in,(R,R)}_{(R,R) \rightarrow (R,P)}$ to the point $(R,P)$ from the point $(R,R)$, to be the set $\{ (z_1,z_2,-1+h,z_4) \ | \ |z_1|, |z_2|, |z_4| <h \}$ 
\end{proof}
\end{lm}

\usetikzlibrary{arrows}
\newcommand{\midarrow}{\tikz \draw[-triangle 90] (0,0) -- +(.1,0);}

\begin{figure}
\centering
\begin{tikzpicture}[xscale=1,>=latex]
       \coordinate  (PP) at (0,0);   
       \coordinate  (RP) at (-4,-2);  \coordinate  (PS) at (4,-2);   
			 \coordinate  (SP) at (4,2); \coordinate (PR) at (-4,2);  
\draw (PP) -- node {\midarrow} (PS);
\draw (RP) -- node {\midarrow} (PP);
\draw (PP) -- node {\midarrow} (SP);
\draw (PR) -- node {\midarrow} (PP);
\node at (0,-0.3) {\scriptsize $(P,P)$};
\node at (-4.2,-2.2) {\scriptsize $(R,P)$};
\node at (4.2,-2.2) {\scriptsize $(P,S)$};
\node at (-4.2,2.2) {\scriptsize $(P,R)$};
\node at (4.2,2.2) {\scriptsize $(S,P)$};
\coordinate (PPSP) at (1,0.5);
\draw (1.4,0.8) -- (1.4,0.1);
\draw[dashed] (1.4,0.1) -- (0.6,0.4);
\draw (0.6,0.4) -- (0.6,1.1);
\draw (0.6,1.1) -- (1.4,0.8);
\node at (1,1.4) {\scriptsize $\Sigma^{out}_{(P,P) \rightarrow (S,P)}$};
\coordinate (PPPS) at (1,-0.5);
\draw[dashed] (1.4,-0.8) -- (1.4,-0.1);
\draw (1.4,-0.1) -- (0.6,-0.4);
\draw (0.6,-1.1) -- (0.6,-0.4);
\draw (0.6,-1.1) -- (1.4,-0.8);
\node at (1,-1.4) {\scriptsize $\Sigma^{out}_{(P,P) \rightarrow (P,S)}$};
\coordinate (PRPP) at (-1,0.5);
\draw (-1.4,0.8) -- (-1.4,0.1);
\draw[dashed] (-1.4,0.1) -- (-0.6,0.4);
\draw (-0.6,0.4) -- (-0.6,1.1);
\draw (-0.6,1.1) -- (-1.4,0.8);
\node at (-1,1.4) {\scriptsize $\Sigma^{in}_{(P,R) \rightarrow (P,P)}$};
\coordinate (RPPP) at (-1,-0.5);
\draw[dashed] (-1.4,-0.8) -- (-1.4,-0.1);
\draw (-1.4,-0.1) -- (-0.6,-0.4);
\draw (-0.6,-1.1) -- (-0.6,-0.4);
\draw (-0.6,-1.1) -- (-1.4,-0.8);
\node at (-1,-1.4) {\scriptsize $\Sigma^{in}_{(R,P) \rightarrow (P,P)}$};
\end{tikzpicture}
\caption{Schematic picture of cross sections near equilibrium $(P,P)$}
\end{figure}
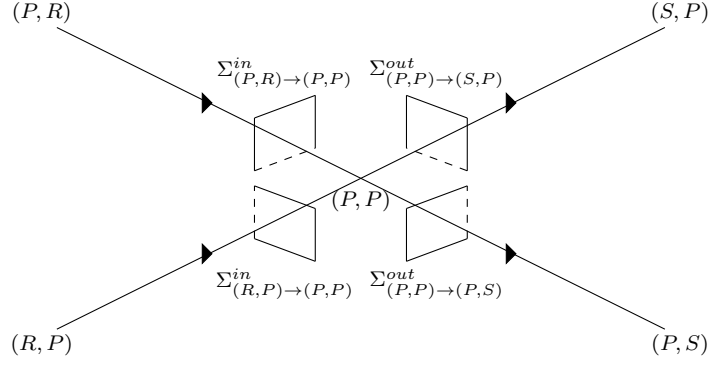

In order to define global and local maps, we need to represent each of the "in cross-sections" in two different coordinates. For the derivation of global maps we will use representation from proposition (\ref{sections}). The following proposition, enables us to compose local maps with global maps along cycles in the network.
\begin{lm}
Representations of the "in" sections in the local coordinates of the points near which they are.
\newline
In $(P,P)$ local coordinates:
\begin{align*}
\Sigma^{in,(P,P)}_{(P,R) \rightarrow (P,P)} = \{ Tran_{(P,R) \rightarrow (P,P)}(z_1,z_2,z_3,z_4)=(z_1, z_2, 1+z_3-z_4, z_4) \ | \ (z_1,z_2,z_3,z_4) \in \Sigma^{in,(P,R)}_{(P,R) \rightarrow (P,P)}  \},
\\
\Sigma^{in,(P,P)}_{(R,P) \rightarrow (P,P)} = \{ Tran_{(R,P) \rightarrow (P,P)}(z_1,z_2,z_3,z_4)=(1+z_1-z_2, z_2, z_3, z_4) \ | \ (z_1,z_2,z_3,z_4) \in \Sigma^{in,(R,P)}_{(R,P) \rightarrow (P,P)}  \},
\end{align*}
In $(P,S)$ local coordinates:
\begin{align*}
\Sigma^{in,(P,S)}_{(P,P) \rightarrow (P,S)} = \{ Tran_{(P,P) \rightarrow (P,S)}(z_1,z_2,z_3,z_4)=(z_1, z_2, z_3, -1+z_3+z_4) \ | \ (z_1,z_2,z_3,z_4) \in \Sigma^{in,(P,P)}_{(P,P) \rightarrow (P,S)}  \},
\\
\Sigma^{in,(P,S)}_{(P,R) \rightarrow (P,S)} = \{ Tran_{(P,R) \rightarrow (P,S)}(z_1,z_2,z_3,z_4)=(z_1, z_2, 1+z_3-z_4,z_3) \ | \ (z_1,z_2,z_3,z_4) \in \Sigma^{in,(P,R)}_{(P,R) \rightarrow (P,S)}  \},
\end{align*}
In $(S,P)$ local coordinates:
\begin{align*}
\Sigma^{in,(S,P)}_{(P,P) \rightarrow (S,P)} = \{ Tran_{(P,P) \rightarrow (S,P)}(z_1,z_2,z_3,z_4)=(z_1, -1+z_1+z_2, z_3, z_4) \ | \ (z_1,z_2,z_3,z_4) \in \Sigma^{in,(P,P)}_{(P,P) \rightarrow (S,P)}  \},
\\ 
\Sigma^{in,(S,P)}_{(R,P) \rightarrow (S,P)} = \{ Tran_{(R,P) \rightarrow (S,P)}(z_1,z_2,z_3,z_4)=(1+z_1-z_2, z_1, z_3, z_4) \ | \ (z_1,z_2,z_3,z_4) \in \Sigma^{in,(R,P)}_{(R,P) \rightarrow (S,P)}  \}
\end{align*}
\begin{proof}
Since the transformation from the $(P,S)$ coordinates to canonical ones is given by the mapping 
$$
z \mapsto MacPS (z+ e_4)
$$ 
and the transformation from canonical to $(S,S)$ is by 
$$z \mapsto MacSS^{-1}(z) - (e_2 +e_4)$$
So the transformation from $(P,S)$ to $(S,S)$ coordinates is given by the mapping 
$$ z \mapsto  -(e_2+e_4) + (MacSS)^{-1} \cdot MacPS (z + e_4)
$$
Hence 
\begin{align*}
\Sigma^{in,(S,S)}_{(P,S) \rightarrow (S,S)} = \{ -(e_2+e_4) + (MacSS)^{-1} \cdot MacPS (z + e_4) \ | \ z \in \Sigma^{in,(P,S)}_{(P,S) \rightarrow (S,S)} \} =
\\ 
= \{ (z_1, -1+z_1+z_2, z_3, z_4) \ | \ (z_1,z_2,z_3,z_4) \in \Sigma^{in,(R,S)}_{(R,S) \rightarrow (R,R)} \}
\end{align*}
\end{proof}
\end{lm}
We find local maps, near each of the equilibrium from the heteroclinic network, in a standard way by linearizing the system near the equilibrium.  
\begin{prop}
Local maps near each point of the heteroclinic network are given by the following formulas: 
\newline
In $(P,S)$ coordinates:
\newline
- from $\Sigma^{in}_{(P,S) \rightarrow (R,S)}$ to $\Sigma^{out}_{(P,S) \rightarrow (R,S)}$, $\Phi_{(P,S) \rightarrow (R,S)} (z_1,z_2,z_3,z_4):= \Big(h, (\frac{h}{z_1})^{\frac{1+ \epsilon_X}{2}}z_2, \frac{z_1 z_3}{h}, (\frac{z_1}{h})^{\frac{1- \epsilon_Y}{2}}z_4 \Big)$
\newline
- from $\Sigma^{in}_{(P,S) \rightarrow (S,S)}$ to $\Sigma^{out}_{(P,S) \rightarrow (S,S)}$, $\Phi_{(P,S) \rightarrow (S,S)} (z_1,z_2,z_3,z_4):= \Big( (\frac{h}{z_2})^{\frac{2}{1+ \epsilon_X}} z_1, h, (\frac{z_2}{h})^{\frac{2}{1+ \epsilon_X}}z_3, (\frac{z_2}{h})^{\frac{1- \epsilon_Y}{1+ \epsilon_X}}z_4 \Big)$
\\ \\
In $(P,P)$ coordinates:
\newline
- from $\Sigma^{in}_{(P,P) \rightarrow (P,S)}$ to $\Sigma^{out}_{(P,P) \rightarrow (P,S)}$, $\Phi_{(P,P) \rightarrow (P,S)} (z_1,z_2,z_3,z_4):= \Big( z_1 (\frac{z_4}{h})^{\frac{1+ \epsilon_X}{1- \epsilon_Y}}, z_2 (\frac{h}{z_4})^{\frac{1- \epsilon_X}{1- \epsilon_Y}}, z_3 (\frac{z_4}{h})^{\frac{1- \epsilon_X}{1- \epsilon_Y}}, h \Big)$
\newline
- from $\Sigma^{in}_{(P,P) \rightarrow (S,P)}$ to $\Sigma^{out}_{(P,P) \rightarrow (S,P)}$, $\Phi_{(P,P) \rightarrow (S,P)} (z_1,z_2,z_3,z_4):= \Big( (\frac{z_2}{h})^{\frac{1+ \epsilon_X}{1- \epsilon_X}} z_1, h, (\frac{z_2}{h})^{\frac{1+ \epsilon_Y}{1- \epsilon_X}} z_3, (\frac{h}{z_2})^{\frac{1- \epsilon_Y}{1- \epsilon_X}}z_4 \Big)$
\\ \\
In $(S,P)$ coordinates:
\newline
- from $\Sigma^{in}_{(S,P) \rightarrow (S,R)}$ to $\Sigma^{out}_{(S,P) \rightarrow (S,R)}$, 
$\Phi_{(S,P) \rightarrow (S,R)} (z_1,z_2,z_3,z_4):= \Big( \frac{z_1 z_3}{h}, (\frac{z_3}{h})^{\frac{1- \epsilon_X}{2}} z_2, h, (\frac{h}{z_3})^{\frac{1+ \epsilon_Y}{2}} z_4 \Big)$
\newline
- from $\Sigma^{in}_{(S,P) \rightarrow (S,S)}$ to $\Sigma^{out}_{(S,P) \rightarrow (S,S)}$, $\Phi_{(S,P) \rightarrow (S,S)} (z_1,z_2,z_3,z_4):= \Big( (\frac{z_4}{h})^{\frac{2}{1+ \epsilon_Y}} z_1, (\frac{z_4}{h})^{\frac{1- \epsilon_X}{1+ \epsilon_Y}} z_2, (\frac{h}{z_4})^{\frac{2}{1+ \epsilon_Y}} z_3, h \Big)$
\begin{proof}
We approximate the system near $(P,S)$ by its linearization, in $(P,S)$ coordinates already: 
$$
\begin{bmatrix}
\dot{z_1} \\
\dot{z_2} \\
\dot{z_3} \\
\dot{z_4}
\end{bmatrix} = 
\begin{bmatrix}
2 && 0 						&& 0 						&& 0 \\
0	&& 1+\epsilon_X && 0 						&& 0 \\
0	&& 0 						&& -2 					&& 0 \\
0 && 0 						&& 0					  && -1+\epsilon_Y
\end{bmatrix} \cdot 
\begin{bmatrix}
z_1 \\
z_2 \\
z_3 \\
z_4
\end{bmatrix} 
$$
The solution to it with initial condition $\big( z_1(0),z_2(0),z_3(0),z_4(0) \big)=(z_1,z_2,z_3,z_4)$ is given by 
$$
\big(z_1(t),z_2(t),z_3(t),z_4(t) \big)= \big( e^{2t}z_1, e^{t(1+\epsilon_X)}z_2, e^{-2t}z_3, e^{t(-1+ \epsilon_Y)}z_4 \big)
$$
For the point $\big( z_1(t),z_2(t),z_3(t),z_4(t) \big)$ to belong to the out section $\Sigma^{out}_{(P,S) \rightarrow (R,S)}$, it is needed to have $z_1(t)=h \Leftrightarrow e^t= (\frac{h}{z_1})^{\frac{1}{2}}$. Hence 
the local map to section $\Sigma^{out}_{(P,S) \rightarrow (R,S)}$ is given by mapping
\begin{multline*}
\Phi_{(P,S) \rightarrow (R,S)} (z_1,z_2,z_3,z_4):= \Big( z_1 ((\frac{h}{z_1})^{\frac{1}{2}})^{2}, z_2 ((\frac{h}{z_1})^{\frac{1}{2}})^{1+ \epsilon_X}, z_3 ((\frac{h}{z_1})^{\frac{1}{2}})^{-2}, z_4 ((\frac{h}{z_1})^{\frac{1}{2}})^{-1+\epsilon_Y} \Big)=
\\
\Big( h, (\frac{h}{z_1})^{\frac{1+ \epsilon_X}{2}}z_2, \frac{z_1 z_3}{h}, (\frac{z_1}{h})^{\frac{1- \epsilon_Y}{2}}z_4 \Big)
\end{multline*}
\end{proof}
\end{prop}
In order to order to introduce global maps from "out cross-sections" to "in cross-sections", we linearise the system along heteroclinic orbits. Note the difference between global maps from the proposition below and those from Aguiar\&Castro paper. Although in section 4 we will find some similarities between the two approaches, we have to point out that we cannot just assume that the global maps are identity as this would lead to the contradiction with mean value theorem - in Aguiar\&Castro paper no orbit from the neighbourhood of the heteroclinic orbit would ever leave it. The actual situation is obviously different, here the global maps derived below are in fact linear - given by the diagonal matrix with some of the entries bigger than 1.  
\begin{prop}
Global maps along each heteroclinic orbit in the network are given by:
\\ \\
In $(P,P)$ coordinates: 
$$\Psi_{(P,P) \rightarrow (P,S)}(z_1, z_2, z_3, h):= \Big( \big(\frac{1-h}{h} \big)^{\frac{1-\epsilon_X}{1-\epsilon_Y}} z_1, \big(\frac{1-h}{h} \big)^{\frac{2}{1-\epsilon_Y}} z_2, \big(\frac{1-h}{h} \big)^{-\frac{3+\epsilon_Y}{1-\epsilon_Y}} z_3,1-h \Big)$$ 
$$\Psi_{(P,P) \rightarrow (S,P)}(z_1, h, z_3, z_4):= \Big( \big(\frac{1-h}{h} \big)^{-\frac{3+\epsilon_X}{1-\epsilon_X}} z_1, 1-h ,\big(\frac{1-h}{h} \big)^{\frac{1-\epsilon_Y}{1-\epsilon_X}} z_3, \big(\frac{1-h}{h} \big)^{-\frac{2}{1-\epsilon_X}} z_4 \Big)$$ 
\newline
In $(P,S)$ coordinates: 
$$\Psi_{(P,S) \rightarrow (R,S)}(h, z_2, z_3, z_4):= \Big(1-h, z_2 \big(\frac{1-h}{h} \big)^{\epsilon_X}, z_3 \big(\frac{1-h}{h} \big)^{- \frac{1-\epsilon_Y}{2}}, z_4 \big(\frac{1-h}{h} \big)^{\frac{1+ \epsilon_Y}{2}}\Big)$$ 
$$\Psi_{(P,S) \rightarrow (S,S)}(z_1, h, z_3, z_4):= \Big(z_1 \big(\frac{1-h}{h} \big)^{- \frac{-3+\epsilon_X}{1+\epsilon_X}}, 1-h, z_3 \big(\frac{1-h}{h} \big)^{- \frac{1+\epsilon_Y}{1+ \epsilon_X}}, z_4 \big(\frac{1-h}{h} \big)^{- \frac{2}{1+ \epsilon_X}} \Big)$$ 
\newline 
In $(S,P)$ coordinates: 
$$\Psi_{(S,P) \rightarrow (S,S)}(z_1, z_2, z_3, h):= \Big(z_1 \big(\frac{1-h}{h} \big)^{- \frac{1+\epsilon_X}{1+\epsilon_Y}}, z_2 \big(\frac{1-h}{h} \big)^{- \frac{2}{1+\epsilon_Y}}, z_3 \big(\frac{1-h}{h} \big)^{-1+ \frac{4}{1+\epsilon_Y}}, 1-h \Big)$$ 
$$\Psi_{(S,P) \rightarrow (S,R)}(z_1, z_2, h, z_4):= \Big( z_1 \big(\frac{1-h}{h} \big)^{- \frac{1- \epsilon_X}{2}}, z_2 \big(\frac{1-h}{h} \big)^{\frac{1+\epsilon_X}{2}}, 1-h, z_4 \big(\frac{1-h}{h} \big)^{\epsilon_Y} \Big)$$ 
\begin{proof}
We compute global maps as shown below in the case of global map $\Psi_{(P,S) \rightarrow (R,S)}$. In the initial equation (\ref{rownanieZredukowane}), let us change the coordinates from the canonical basis to local $(P,S)$ coordinates $(z_1,z_2,z_3,z_4)$ in which $(P,S)$ is an origin, and $(R,S)_{(P,S)}=(1,0,0,0)$. Let us rewrite the right-hand side of the system (after changing the coordinates) as 
$$
G(z_1,z_2,z_3,z_4)= \big( \zeta (z_1,z_2,z_3,z_4), g(z_1,z_2,z_3,z_4) \big) \in \mathbb R \times \mathbb R^3
$$ 
and the heteroclinic orbit from $(P,S)$ to $(R,S)$ by $\alpha= \{ (z_1,0,0,0) \ | \ z_1 \in (0,1) \}$. Near $\alpha$ one can write $(z_1,z_2,z_3,z_4)=(z_1,z)$ with $z \in \mathbb R^3$ -small. Then 
$$
z'(t) = g(z_1(t),z(t))= g(z_1(t),0)+g_z(z_1(t),0).z(t) + h.o.t.
$$ 
Since $g(z_1(t),0)=0$ (because of the invariance of $\alpha$) and by the chain rule 
$$
z'(t)= \frac{\partial z}{\partial z_1} \cdot z_1 '(t)
$$
we have that 
\begin{multline}\label{linAlongOrbit}
\frac{\partial z}{\partial z_1} \approx \frac{1}{z_1'} g_z(z_1,0).z = \frac{1}{\zeta (z_1,0,0,0)} g_z(z_1,0).z =
\\
= \frac{1}{2(1-z_1)z_1} \cdot 
\begin{bmatrix}
1+\epsilon_X -2z_1 && 0 										&& 0 \\
0 						     && -2+(3+\epsilon_Y) z_1 && 0 \\
0 						     && 0					  					&& -1+\epsilon_Y (1-z_1)+3z_1
\end{bmatrix} \cdot 
\begin{bmatrix}
z_2 \\
z_3 \\
z_4
\end{bmatrix}  
\end{multline}
\newline
Imposing initial condition $z(h)=z_0$, 
$$
z(\tau)= Exp \Big( \int_{h}^{\tau}{\frac{1}{\zeta (z_1,0)} g_z(z_1,0) dz_1} \Big).z_0
$$ 
is the solution to the latter equation (\ref{linAlongOrbit}). The global map $\Psi_{(P,S) \rightarrow (R,S)}$ will be given by mapping 
$$
\big(h,z(h) \big) \longmapsto \big(1-h,z(1-h) \big)
$$ 
\newline 
Hence in $(P,S)$ coordinates global map to $(R,S)$ has the form: 
$$
\Psi_{(P,S) \rightarrow (R,S)}(h, z_2, z_3, z_4):= \Big( 1-h, z_2 \big(\frac{1-h}{h} \big)^{\epsilon_X}, z_3 \big(\frac{1-h}{h} \big)^{-\frac{1-\epsilon_Y}{2}}, z_4 \big(\frac{1-h}{h} \big)^{\frac{1+\epsilon_Y}{2}} \Big)
$$ 
with $z_1=h$ fixed.
\end{proof}
\end{prop}

\begin{prop}
The return map from any section within invariant square $\big( (R,P),(S,P),(S,S),(R,S) \big)$ to itself is an identity map (but restricted to a smaller neighbourhood of origin in the cross section).
The same observation holds for the return maps within other invariant squares: 
\newline 
$ \big( (R,P),(S,P),(S,R),(R,R) \big)$, $\big( (P,R),(P,S),(S,S),(S,R) \big)$, $\big( (P,P),(S,P),(S,R),(P,R) \big)$, 
\newline
$\big( (P,P),(P,S),(R,S),(R,P) \big)$, $\big( (R,R),(P,R),(P,S),(R,S) \big)$
\begin{proof}
It suffices to show that the return map to the cross section $\Sigma^{out}_{(R,P) \rightarrow (S,P)}$  within invariant square $((R,P),(S,P),(S,S),(R,S))$ is an identity. For this, we compute
\begin{multline*} 
\Phi_{(R,P) \rightarrow (S,P)} \circ Tran_{(R,S) \rightarrow (R,P)} \circ \Psi_{(R,S) \rightarrow (R,P)} \circ \Phi_{(R,S) \rightarrow (R,P)} \circ Tran_{(S,S) \rightarrow (R,S)} \circ \Psi_{(S,S) \rightarrow (R,S)} \circ \Phi_{(S,S) \rightarrow (R,S)} \circ 
\\
\circ Tran_{(S,P) \rightarrow (S,S)} \circ \Psi_{(S,P) \rightarrow (S,S)} \circ \Phi_{(S,P) \rightarrow (S,S)} \circ Tran_{(R,P) \rightarrow (S,P)} \circ \Psi_{(R,P) \rightarrow (S,P)}(0, h, 0, z_4)=
\end{multline*}
$$
\Big ( 0,h,0, (\frac{1-h}{h})^{\frac{1+\epsilon_X}{2}} h \Big ( (\frac{1-h}{h})^{\frac{-1+\epsilon_Y}{-1+\epsilon_X}} \big ( (\frac{1-h}{h})^{-\frac{1+\epsilon_X}{1+\epsilon_Y}} ((\frac{h}{1-h})^{ \frac{1- \epsilon_Y}{2}} \frac{z_4}{h}))^{\frac{2}{1+\epsilon_Y}} \big )^{\frac{1+\epsilon_Y}{1-\epsilon_X}} \Big )^{\frac{1-\epsilon_X}{2}} \Big ) =
$$
$$
\Big( 0,h,0, (\frac{1-h}{h})^{\frac{1+\epsilon_X}{2}+\frac{1-\epsilon_Y}{1-\epsilon_X} \frac{1-\epsilon_X}{2} - \frac{1+\epsilon_X}{1+\epsilon_Y} \frac{1+\epsilon_Y}{1-\epsilon_X} \frac{1-\epsilon_X}{2} - \frac{1-\epsilon_Y}{2} \frac{2}{1+\epsilon_Y} \frac{1+\epsilon_Y}{1-\epsilon_X} \frac{1-\epsilon_X}{2}} \cdot h^{1 - \frac{2}{1+\epsilon_Y} \frac{1+\epsilon_Y}{1-\epsilon_X} \frac{1-\epsilon_X}{2}} \cdot z_4 \Big)=(0,h,0,z_4)
$$ 
\end{proof}
\end{prop}

\section{Asymptotic dynamics near the network}
In this section we will provide analytical proofs of numerical findings by Sato et al. \cite{[17]} and proper analysis of the behaviour close to the heteroclinic network. In particular we will explain the differences and similarities between our model and the one built by Aguiar\&Castro.
\\ \\
Recall that:
$$
C_0 \ cycle = PS \rightarrow RS \rightarrow RP \rightarrow SP \rightarrow SR \rightarrow PR \rightarrow PS
$$
$$
C_1 \ cycle = PP \rightarrow SP \rightarrow SS \rightarrow RS \rightarrow RR \rightarrow PR \rightarrow PP
$$
$$
C_2 \ cycle = PP \rightarrow PS \rightarrow SS \rightarrow SR \rightarrow RR \rightarrow RP \rightarrow PP
$$

\subsection{Dynamics near $C_0$, $C_1$, $C_2$ cycles}
We start with an easy observation, that symmetry appearing in the system simplifies our computations along $C_0$, $C_1$, $C_2$ cycles, so instead of composing maps along six heteroclinic orbits, we can investigate part of it along two heteroclinic orbits.
\begin{prop}\label{returnMaps}
The return maps to the sections $\Sigma^{out}_{(P,P) \rightarrow (S,P)}$, $\Sigma^{out}_{(P,P) \rightarrow (P,S)}$, $\Sigma^{out}_{(P,S) \rightarrow (R,S)}$ through the $C_1$, $C_2$, $C_0$  cycles, respectively, are given by:
\\ \\
$Ret_{(P,P), C_1}= (\sigma^{-1} \circ Tran_{(S,S) \rightarrow (P,P)} \circ \Phi_{(S,S) \rightarrow (R,S)} \circ $
\begin{equation}\label{returnMapC1}
\circ \ Tran_{(S,P) \rightarrow (S,S)} \circ \Psi_{(S,P) \rightarrow (S,S)} \circ
\Phi_{(S,P) \rightarrow (S,S)} \circ Tran_{(P,P) \rightarrow (S,P)} \circ \Psi_{(P,P) \rightarrow (S,P)})^3
\end{equation}
$Ret_{(P,P), C_2}= (\sigma^{-1} \circ Tran_{(S,S) \rightarrow (P,P)} \circ \Phi_{(S,S) \rightarrow (S,R)} \circ $
\begin{equation}\label{returnMapC2}
\circ \ Tran_{(P,S) \rightarrow (S,S)} \circ \Psi_{(P,S) \rightarrow (S,S)} \circ \Phi_{(P,S) \rightarrow (S,S)} \circ Tran_{(P,P) \rightarrow (P,S)} \circ \Psi_{(P,P) \rightarrow (P,S)})^3
\end{equation}
$Ret_{(P,S), C_0}= (Tran_{(P,P) \rightarrow (P,S)} \circ \sigma \circ Tran_{(R,P) \rightarrow (P,P)} \circ \Phi_{(R,P) \rightarrow (S,P)} \circ $
\begin{equation}\label{returnMapC0}
\circ \ Tran_{(R,S) \rightarrow (R,P)} \circ \Psi_{(R,S) \rightarrow (R,P)} \circ \Phi_{(R,S) \rightarrow (R,P)} \circ Tran_{(P,S) \rightarrow (R,S)} \circ \Psi_{(P,S) \rightarrow (R,S)})^3
\end{equation}
\begin{proof}
For $C_2$ cycle, let 
$$
T_{(P,P) \rightarrow (S,S)}:=\Phi_{(R,P) \rightarrow (S,P)} \circ Tran_{(R,S) \rightarrow (R,P)} \circ \Psi_{(R,S) \rightarrow (R,P)} \circ \Phi_{(R,S) \rightarrow (R,P)} \circ Tran_{(P,S) \rightarrow (R,S)} \circ \Psi_{(P,S) \rightarrow (R,S)}
$$
then 
$$
Ret_{(P,P), C_2}= T_{(R,R) \rightarrow (P,P)} \circ T_{(S,S) \rightarrow (R,R)} \circ T_{(P,P) \rightarrow (S,S)}
$$
with $T_{(R,R) \rightarrow (P,P)}$, $T_{(S,S) \rightarrow (R,R)}$ defined analogically (through $C_2$ cycle).
Due to the symmetry given by the action $\sigma$ one can write: 
\newline
$$Ret_{(P,P), C_2}= \big( T_{(R,R) \rightarrow (P,P)} \circ Tran_{(P,P) \rightarrow (R,R)} \big) \circ \big(Tran_{(R,R) \rightarrow (P,P)}  \circ T_{(S,S) \rightarrow (R,R)} \circ Tran_{(P,P) \rightarrow (S,S)} \big) \ \circ 
$$
$$
\circ \ \big(Tran_{(S,S) \rightarrow (P,P)}  \circ T_{(P,P) \rightarrow (S,S)} \big) = 
\big(\sigma^2 \circ Tran_{(S,S) \rightarrow (P,P)} \circ T_{(P,P) \rightarrow (S,S)} \circ \sigma^{-2} \big) \ \circ
$$
$$
\circ \ \big(\sigma \circ Tran_{(S,S) \rightarrow (P,P)} \circ T_{(P,P) \rightarrow (S,S)} \circ \sigma^{-1} \big) \circ \big( Tran_{(S,S) \rightarrow (P,P)} \circ T_{(P,P) \rightarrow (S,S)} \big) =$$ 
$$
= \big( \sigma^{-1} \circ Tran_{(S,S) \rightarrow (P,P)} \circ T_{(P,P) \rightarrow (S,S)} \big)^3
$$
\end{proof}
\end{prop}
Let us identify, by symmetry $\sigma$, tie:=$\{ (P,P), (S,S), (R,R) \} / \sigma$, win/loss:=$\{ (S,P), (P,R), (R,S) \} / \sigma$, loss/win:=$\{ (P,S), (S,R), (R,P) \} / \sigma$. Hence we can enumerate "out cross-sections" as in the figure 4. 
Now, if $h<<1$ is set to be a constant we can find approximations of the maps between the sections 1-6 in the limit of coordinates on the cross-sections tending to zero. 

\begin{prop}\label{CAFormulas}
With $h<<1$ and $0 < s_1, s_2, s_3, s_4 <\epsilon <<h$, it is possible to rescale coordinates linearly at the cross sections near: $tie$, $loss/win$, $win/loss$, so that the mappings are of the "constant-free" form 
\newline
$\vartheta_{12}(s_1,s_2,s_3,h)= \big( s_1 s_2^{-\frac{2}{1+\epsilon_X}},h,s_2^{\frac{2}{1+\epsilon_X}} s_3, s_2^{\frac{1-\epsilon_Y}{1+\epsilon_X}} \big) = $
$$
=\Phi_{(P,S) \rightarrow (S,S)} \circ Tran_{(P,P) \rightarrow (P,S)} \circ \Phi_{(P,P) \rightarrow (P,S)}(s_1,s_2,s_3,h)
$$ 
$\vartheta_{14}(s_1,s_2,s_3,h)= \big( h, s_1^{-\frac{1+\epsilon_X}{2}} s_2, s_1 s_3, s_1^{\frac{1-\epsilon_Y}{2}} \big)=$
$$
= \Phi_{(P,S) \rightarrow (R,S)} \circ Tran_{(P,P) \rightarrow (P,S)} \circ \Phi_{(P,P) \rightarrow (P,S)}(s_1,s_2,s_3,h)
$$
$\vartheta_{21}(s_1,h,s_3,s_4)= \big( s_3^{\frac{1+\epsilon_X}{1-\epsilon_Y}}, s_1 s_3^{-\frac{1-\epsilon_X}{1-\epsilon_Y}}, s_3^{\frac{1+\epsilon_Y}{1-\epsilon_Y}} s_4, h \big)= $
$$ 
=\sigma^{-1} \circ Tran_{(S,S) \rightarrow (P,P)} \circ \Phi_{(S,S) \rightarrow (S,R)} \circ Tran_{(P,S) \rightarrow (S,S)} \circ \Phi_{(P,S) \rightarrow (S,S)} (s_1,h,s_3,s_4)
$$
$\vartheta_{26}(s_1,h,s_3,s_4)= \big( s_1^{\frac{1+\epsilon_X}{1-\epsilon_X}},h, s_1^{\frac{1+\epsilon_Y}{1-\epsilon_X}} s_4, s_1^{-\frac{1-\epsilon_Y}{1-\epsilon_X}} s_3 \big) = $
$$
= \sigma^{-1} \circ Tran_{(S,S) \rightarrow (P,P)} \circ  \Phi_{(S,S) \rightarrow (R,S)} \circ Tran_{(P,S) \rightarrow (S,S)} \circ \Phi_{(P,S) \rightarrow (S,S)}(s_1,h,s_3,s_4)
$$
$\vartheta_{32}(s_1,s_2,h,s_4)= \big( s_1^{-\frac{2}{1+\epsilon_X}} s_2, h, s_1^{\frac{2}{1+\epsilon_X}}, s_1^{\frac{1-\epsilon_Y}{1-\epsilon_X}} s_4 \big)= $
$$
= Tran_{(P,P) \rightarrow (P,S)} \circ \sigma^{-1} \circ Tran_{(S,R) \rightarrow (P,P)} \circ \Phi_{(S,R) \rightarrow (R,R)} \circ Tran_{(S,P) \rightarrow (S,R)} \circ \Phi_{(S,P) \rightarrow (S,R)}(s_1,s_2,h,s_4)
$$
$\vartheta_{34}(s_1,s_2,h,s_4)= \big( h, s_1 s_2^{-\frac{1+\epsilon_X}{2}}, s_2, s_2^{\frac{1-\epsilon_Y}{2}} s_4 \big)= $
$$
= Tran_{(P,P) \rightarrow (P,S)} \circ \sigma^{-1} \circ Tran_{(S,R) \rightarrow (P,P)} \circ \Phi_{(S,R) \rightarrow (P,R)} \circ Tran_{(S,P) \rightarrow (S,R)} \circ \Phi_{(S,P) \rightarrow (S,R)}(s_1,s_2,h,s_4)
$$
$\vartheta_{43}(h,s_2,s_3,s_4)= \big( s_4, s_2 s_4^{\frac{1-\epsilon_X}{2}}, h, s_3 s_4^{-\frac{1+\epsilon_Y}{2}} \big)=$
$$
= Tran_{(P,P) \rightarrow (S,P)} \circ \sigma^{-1} \circ Tran_{(R,S) \rightarrow (P,P)} \circ \Phi_{(R,S) \rightarrow (R,P)} \circ Tran_{(P,S) \rightarrow (R,S)} \circ \Phi_{(P,S) \rightarrow (R,S)}(h,s_2,s_3,s_4)
$$
$\vartheta_{45}(h,s_2,s_3,s_4)= \big( s_3^{\frac{2}{1+\epsilon_Y}}, s_2 s_3^{\frac{1-\epsilon_X}{1+\epsilon_Y}} , s_3^{-\frac{2}{1+\epsilon_Y}} s_4, h \big)= $
$$
= Tran_{(P,P) \rightarrow (S,P)} \circ \sigma^{-1} \circ Tran_{(S,R) \rightarrow (P,P)} \circ \Phi_{(R,S) \rightarrow (R,R)} \circ Tran_{(P,S) \rightarrow (R,S)} \circ \Phi_{(P,S) \rightarrow (R,S)}(h,s_2,s_3,s_4)
$$
$\vartheta_{51}(s_1,s_2,s_3,h)= \big( s_2 s_3^{\frac{1+\epsilon_X}{1-\epsilon_Y}}, s_1 s_3^{-\frac{1-\epsilon_X}{1-\epsilon_Y}} , s_3^{\frac{1+\epsilon_Y}{1-\epsilon_Y}}, h \big)= $
$$= \sigma^{-1} \circ Tran_{(S,S) \rightarrow (P,P)} \circ \Phi_{(S,S) \rightarrow (S,R)} \circ Tran_{(S,P) \rightarrow (S,S)} \circ \Phi_{(S,P) \rightarrow (S,S)}(s_1,s_2,s_3,h)
$$
$\vartheta_{56}(s_1,s_2,s_3,h)= \big( s_2 s_3^{\frac{1+\epsilon_X}{1-\epsilon_Y}}, s_1 s_3^{-\frac{1-\epsilon_X}{1-\epsilon_Y}} , s_3^{\frac{1+\epsilon_Y}{1-\epsilon_Y}}, h \big) = $
$$
= \sigma^{-1} \circ Tran_{(S,S) \rightarrow (P,P)} \circ \Phi_{(S,S) \rightarrow (R,S)} \circ Tran_{(S,P) \rightarrow (S,S)} \circ \Phi_{(S,P) \rightarrow (S,S)}(s_1,s_2,s_3,h)
$$
$\vartheta_{63}(s_1,h,s_3,s_4)= \big( s_1 s_3, s_3^{\frac{1-\epsilon_X}{2}}, h, s_3^{-\frac{1+\epsilon_Y}{2}} s_4 \big)= $
$$\Phi_{(S,P) \rightarrow (S,R)} \circ Tran_{(P,P) \rightarrow (S,P)} \circ \Phi_{(P,P) \rightarrow (S,P)}(s_1,h,s_3,s_4)
$$
$\vartheta_{65}(s_1,h,s_3,s_4)= \big( s_1 s_4^{\frac{2}{1+\epsilon_Y}}, s_4^{\frac{1-\epsilon_X}{1+\epsilon_Y}}, s_3 s_4^{-\frac{2}{1+\epsilon_Y}}, h \big)= $
$$=\Phi_{(S,P) \rightarrow (S,S)} \circ Tran_{(P,P) \rightarrow (S,P)} \circ \Phi_{(P,P) \rightarrow (S,P)}(s_1,h,s_3,s_4)
$$
\end{prop}

\begin{figure}\label{QN}
\centering
\begin{tikzpicture}[xscale=1,>=latex]
\coordinate (WL) at (-4,0);
\coordinate (LW) at (4,0);
\coordinate (D) at (0,6.92820324);
\filldraw
(WL) circle (2pt);
\filldraw
(LW) circle (2pt);
\filldraw
(D) circle (2pt);
\node at (-4.5,0) {\scriptsize $W/L$};
\node at (4.5,0) {\scriptsize $L/W$};
\node at (0,7.2) {\scriptsize $D$};
\draw[decoration={markings, mark=at position 0.5 with {\arrow{>}}}, postaction={decorate}] (WL) to[out=-25,in=205] (LW);
\draw[decoration={markings, mark=at position 0.5 with {\arrow{>}}}, postaction={decorate}] (LW) to[out=170,in=10] (WL);
\draw[decoration={markings, mark=at position 0.5 with {\arrow{>}}}, postaction={decorate}] (WL) to[out=50,in=250] (D);
\draw[decoration={markings, mark=at position 0.5 with {\arrow{>}}}, postaction={decorate}] (D) to[out=210,in=90] (WL);
\draw[decoration={markings, mark=at position 0.5 with {\arrow{>}}}, postaction={decorate}] (LW) to[out=90,in=-30] (D);
\draw[decoration={markings, mark=at position 0.5 with {\arrow{>}}}, postaction={decorate}] (D) to[out=-70,in=130] (LW);
\draw (3.6,1) -- (4.5,1);
\draw[dashed] (3.4,1.5) -- (4.3,1.5);
\draw (3.6,1) -- (3.4,1.5);
\draw (4.3,1.5) -- (4.5,1);
\node at (4.7,1.3) {\scriptsize $2$};
\draw (0,6) -- (0.8,6);
\draw[dashed] (0.2,5.5) -- (1,5.5);
\draw (0,6) -- (0.2,5.5);
\draw (0.8,6) -- (1,5.5);
\node at (0.1,5.3) {\scriptsize $1$};
\draw (-1.8,6) -- (-1,6);
\draw[dashed] (-2.2,5.5) -- (-1.4,5.5);
\draw (-1.8,6) -- (-2.2,5.5);
\draw (-1,6) -- (-1.4,5.5);
\node at (-2.4,5.7) {\scriptsize $6$};
\draw (2.9,-0.2) -- (2.9,0.6);
\draw[dashed] (2.4,0) -- (2.4,0.8);
\draw (2.9,-0.2) -- (2.4,0);
\draw (2.9,0.6) -- (2.4,0.8);
\node at (2.7,0.9) {\scriptsize $4$};
\draw (-3,-0.1) -- (-3,-0.8);
\draw[dashed] (-2.5,-0.3) -- (-2.5,-1);
\draw (-3,-0.1) -- (-2.5,-0.3);
\draw (-3,-0.8) -- (-2.5,-1);
\node at (-2.3,-0.2) {\scriptsize $3$};
\draw (-3.5,1) -- (-2.7,1);
\draw[dashed] (-3.3,1.5) -- (-2.5,1.5);
\draw (-3.5,1) -- (-3.3,1.5);
\draw (-2.7,1) -- (-2.5,1.5);
\node at (-2.4,1.25) {\scriptsize $5$};
\end{tikzpicture}
\caption{Schematic picture of the quotient network together with cross-sections (1-6) after identification of equilibria by $\sigma$}
\end{figure}
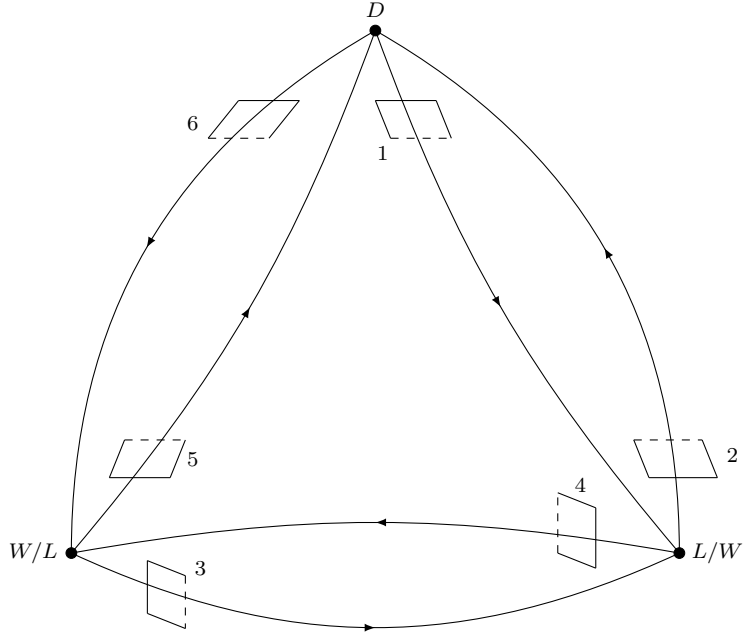
\begin{rmq}
Note that in the proposition above, we have obtained the same formulas as Aguiar\&Castro (upto some permutation, that comes from the difference between our and their transition maps), however the major difference is in the form of the domains of these maps (precisely: the constants appearing in the formulas), which in fact we will not use in the later analysis. 
\end{rmq}
Before stating next proposition, we need to recall some definitions from \cite{[12]} and \cite{[13]}.

\begin{df}
A compact invariant set $X$ is called:
\newline
- asymptotically stable relative to $N$ if for every neighbourhood $U$ of $X$ there is a neighbourhood $V$ of $X$ such that for all $x \in V \cap N$ we have $\omega(x) \subset X$ and $\phi_t(x) \in U \cap N$ for all $t>0$.  
\newline 
- completely unstable relative to $N$ if there is a neighbourhood $U$ of $X$ such that for all $x \in U \cap N$ there is $t_0>0$ with $\phi_{t_0}(x) \notin U$.
\end{df}
\begin{df}
A compact invariant set $X \subset \mathbb R^n$ is called essentially asymptotically stable $(e.a.s.)$ if it is asymptotically stable relative to a set $N \subset \mathbb R^n$ with the property that: 
\begin{equation}\label{relativeStability}
\lim_{\epsilon \rightarrow 0} \frac{ \mathcal{L}^n (\{ x \in \mathbb R^n \ | \ dist(x,X) <\epsilon \}  \cap N)}{ \mathcal{L}^n (\{ x \in \mathbb R^n \ | \ dist(x,X) <\epsilon \} )} =1
\end{equation}  
\end{df}
\begin{df}
A compact invariant set $X$ is almost completely unstable (a.c.u.) if it is completely unstable relative to a set $N \subset \mathbb R^n$ with property (\ref{relativeStability}).
\end{df}
\begin{df}
For $x \in X$, and small enough $\epsilon, \delta >0$ define
$$
\Sigma_{\epsilon,\delta}(x):= \frac{ \mathcal{L}^n \big( B(x,\epsilon)  \cap \mathcal{B_\delta}(X) \big) }{\mathcal{L}^n \big( B(x,\epsilon) \big)}
$$  
where $\mathcal{B_\delta}(X) = \Big{\{ } x \in \mathcal{B}(X) \ | \ \forall t>0: \phi_t(x) \in \{ z \in \mathbb R^n \ | \ dist(z,X) <\delta \} \Big { \} }$  and $\mathcal{B}(X)$ is the basin of attraction of $X$, i.e. the set of points $x \in \mathbb R^n$ with $\omega(x) \subset X$.
\end{df}
\begin{df}
Local stability index at $x \in X$ is defined to be:
$$
\sigma_{loc}(x):=\sigma_{loc,+}(x) - \sigma_{loc,-}(x)
$$
where
$$
\sigma_{loc,-}(x):= \lim_{\delta \rightarrow 0} \lim_{\epsilon \rightarrow 0} \ \frac{ln \big(\Sigma_{\epsilon,\delta} (x) \big)}{ln(\epsilon)} \ \ \ \ \ \
\sigma_{loc,+}(x):= \lim_{\delta \rightarrow 0} \lim_{\epsilon \rightarrow 0} \ \frac{ln \big(1-\Sigma_{\epsilon,\delta} (x) \big)}{ln(\epsilon)}
$$
\end{df}
We will make use of the following theorem (Theorem 3.5 in \cite{[12]}), in order to prove, for some parameter values $\epsilon_X, \epsilon_Y$, essential asymptotic stability of $C_0$ and relative asymptotic stability of $C_1$, $C_2$.
\begin{theo}\label{ALohse} (A.Lohse)
\newline
Let $X \subset \mathbb R^n$ be a heteroclinic cycle or network with finitely many equilibria and connecting trajectories. Suppose that $\mathcal{L}^1 (X)< +\infty$ and that the local stability index $\sigma_{loc}(x)$ exists and is not equal to zero for all $x \in X$. Then generically we have:
\newline
(i) $X$ is e.a.s. $\Leftrightarrow$ $\sigma_{loc}>0$ along all connecting trajectories
\newline
(i) $X$ is a.c.u. $\Leftrightarrow$ $\sigma_{loc}<0$ along all connecting trajectories
\end{theo}
We prove here corrected statement of the first part of the Theorem 4.8 from Aguiar\&Castro paper:  
\begin{thm}\label{C0eas}
For $\epsilon_X+\epsilon_Y<0$, the cycle $C_0$ is essentially asymptotically stable.
\begin{proof}
To prove that $C_0$ is e.a.s., one has to show that all of the local indices along the $C_0$ cycle are positive. From propositions \ref{returnMaps} and \ref{CAFormulas}, in the limit with $h<<1$ and $s_1,s_2,s_3,s_4 < \epsilon <<h$, the formulas for the maps along $C_0$ cycle simplify to the formulas from the paper by Aguiar\&Castro. Although the domains of these maps are different by constants from the domains of maps in Aguilar\&Castro paper, however for the computation of indices the only thing that counts, in this example, is the relationship between the exponents of coordinates in the maps and their domains as well - and this remains the same while taking $s_1,s_2,s_3,s_4 < \epsilon <<h$ and passing to the limit $\epsilon \rightarrow 0$.
\newline
Furthermore, note that in the presence of symmetry $\sigma$ in the system and formula (\ref{returnMapC0}), it suffices to compute only 2 indices for the map 
$$
Tran_{(P,P) \rightarrow (P,S)} \circ \sigma \circ Tran_{(R,P) \rightarrow (P,P)} \circ \Phi_{(R,P) \rightarrow (S,P)} \circ Tran_{(R,S) \rightarrow (R,P)} \circ \Psi_{(R,S) \rightarrow (R,P)} \circ 
$$
$$
\circ \  \Phi_{(R,S) \rightarrow (R,P)} \circ Tran_{(P,S) \rightarrow (R,S)} \circ \Psi_{(P,S) \rightarrow (R,S)}
$$ 
instead of 6 indices for the return map $Ret_{(P,S), C_0}$.
\newline 
Liliana Garrido Silva \cite{[8]} was first that has computed these local indices for the Aguiar\&Castro model and they are all positive for $\epsilon_X+\epsilon_Y<0$. Hence from the above reasoning, and by theorem (\ref{ALohse}), $C_0$ is an e.a.s cycle.
\end{proof}
\end{thm}

Because of the Proposition \ref{CAFormulas}, we can introduce logarithmic coordinates, useful for description of stability of the cycles and bifurcations happening in the system with $\epsilon_X, \epsilon_Y$ varying.

\begin{prop}\label{C0eigen}
The map from $(P,S)$ to $(P,S)$ through $(R,S)$ (and after identifying $(R,P)$ with $(P,S)$ by symmetry $\sigma$) : 
$$
Tran_{(P,P) \rightarrow (P,S)} \circ \sigma \circ Tran_{(R,P) \rightarrow (P,P)} \circ T_{(P,S) \rightarrow (R,S) \rightarrow (R,P)}$$
in logarithmic coordinates $w_i:=- \ln(s_i)$, $i=2,3,4$ has the form:
$$(w_2,w_3,w_4)^T=A_{C_0} \cdot (w_2,w_3,w_4)^T$$
where
$$A_{C_0} = \begin{bmatrix}
-\frac{1+\epsilon_X}{2} && 0 && \frac{3+\epsilon_X^2}{4}  \\
1 && 0 && \frac{1-\epsilon_X}{2}  \\
\frac{1-\epsilon_Y}{2} && 1 && -\frac{1+\epsilon_Y}{2} + \frac{(1-\epsilon_X)(1-\epsilon_Y)}{4} \\
\end{bmatrix}$$
\newline
The eigenvalues $\lambda_i$ and corresponding eigenvectors $v_i=(v_{i1},v_{i2},v_{i3})$, $i=1,2,3$, of matrix $A_{C_0}$ satisfy the following conditions:
\newline 
- $\lambda_1>0$, while $\lambda_2,\lambda_3$ are complex conjugated and $Re(\lambda_2)<0$, $Re(\lambda_3)<0$, $Im(\lambda_2)<0$, $Im(\lambda_3)>0$
\newline
-  $v_{13}=v_{23}=v_{33}=1$ 
\newline
-  $v_{11}>0, v_{12}>0$
\newline
-  $Re(v_{21})<0, Re(v_{22})<0$, $Im(v_{21})>0$, $Im(v_{22})<0$  
\newline
-  $Re(v_{31})<0, Re(v_{32})<0$, $Im(v_{31})<0$, $Im(v_{32})>0$
\newline
-  for $\epsilon_X + \epsilon_Y<0$: 
$$
\lambda_1 > 1 > Abs(\lambda_2)=Abs(\lambda_3)
$$ 
-  for $\epsilon_X + \epsilon_Y > 0$: 
$$
\lambda_1 < 1 < Abs(\lambda_2)=Abs(\lambda_3)
$$
\end{prop}


\begin{figure}
\centering
\includegraphics[width=50mm]{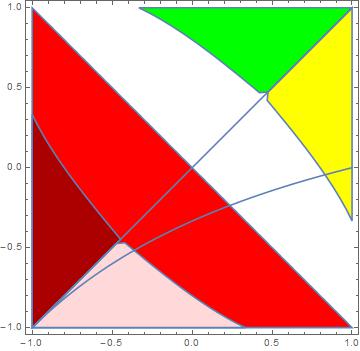}
\caption{The picture shows regions $O_1$, $O_2$, $\epsilon_X+\epsilon_Y<0$, and $\beta$ -curve $\epsilon_X = \frac{1+ 3 \epsilon_Y}{1- \epsilon_Y}$ (where $\beta=\frac{1+ \epsilon_X}{2}$). The $\beta$-curve intersects boundary of the yellow region $O_2$ at point $\beta = \tau \approx 0.915$.  
}
\end{figure}

\begin{prop}\label{C2eigen}
The map from $(P,P)$ to $(P,P)$ through $(P,S)$ (and after identifying $(S,S)$ with $(P,P)$ by symmetry $\sigma$): 
$$\sigma^{-1} \circ Tran_{(S,S) \rightarrow (P,P)} \circ T_{(P,P) \rightarrow (P,S) \rightarrow (S,S)}
$$ 
in logarithmic coordinates 
$w_i:=- \ln(s_i)$, $i=1,2,3$ is given by:
$$
(w_1,w_2,w_3)^T=A_{C_2} \cdot (w_1,w_2,w_3)^T
$$
where 
$$A_{C_2} = \begin{bmatrix}
0 && \frac{2}{1-\epsilon_Y} && \frac{1+\epsilon_X}{1-\epsilon_Y}  \\
1 && - \frac{2}{1+\epsilon_X} - \frac{2 (1-\epsilon_X)}{(1+\epsilon_X)(1-\epsilon_Y)} && -\frac{(1-\epsilon_X)}{(1-\epsilon_Y)}  \\
0 && \frac{(1-\epsilon_Y)}{(1+\epsilon_X)} + \frac{2 (1+\epsilon_Y)}{(1+\epsilon_X)(1-\epsilon_Y)} && \frac{(1+\epsilon_Y)}{(1-\epsilon_Y)}
\end{bmatrix}
$$
The eigenvalues $\lambda_i$ and corresponding eigenvectors $v_i=(v_{i1},v_{i2},v_{i3})$, $i=1,2,3$, of matrix $A_{C_2}$ satisfy the following conditions:
\newline
-  $v_{13}=v_{23}=v_{33}=1$ 
\newline
-  $v_{21}<0, v_{22}<0$
\newline
-  $v_{31}>0, v_{32}<0$
\newline
-  for $\epsilon_X>\epsilon_Y$:  
$$v_{11}>0, \ \ v_{12}>0$$
-  for $\epsilon_X<\epsilon_Y$:  
$$v_{11}>0, \ \ v_{12}<0$$
-  for $\epsilon_X>\epsilon_Y$:    
$$\lambda_2 \in (-1,0), \ \ \ \lambda_3 <-1, \ \ \ \lambda_1>0$$
-  in the yellow region
$$
\lambda_1>max(1, |\lambda_2|, |\lambda_3|)
$$
-  in the white region where $\epsilon_X > \epsilon_Y$
$$
1<\lambda_1<|\lambda_3|
$$
-  in the red region where $\epsilon_X > \epsilon_Y$
$$
\lambda_1<1
$$
- in the light red region where $\epsilon_X > \epsilon_Y$
$$
\lambda_1 < | \lambda_2 |
$$
\end{prop}

\begin{rmq}
We obtain the same proposition, as Prop. \ref{C2eigen}, for cycle $C_1$ by swapping $\epsilon_X$ with $\epsilon_Y$.
\end{rmq}
We are ready to prove that $C_2$ and $C_1$ are asymptotically stable relatively to some open set, however they cannot be essentially asymptotically stable according to the Theorem \ref{ALohse}.
\begin{thm}\label{C2attract}
For some open regions $O_1$ (green), $O_2$ (yellow) in parameter space $\epsilon_X, \epsilon_Y$, cycles $C_1$ and $C_2$, respectively, are neither almost completely unstable nor essentially asymptotically stable. For each of them, there exist an open set of points in the phase space which is attracted to the respective cycle.
\begin{proof}
By propositions \ref{returnMaps} and \ref{C2eigen}, the return map to the cross section $\Sigma^{out}_{(P,P) \rightarrow (P,S)}$ through the whole $C_2$ cycle is given by 
$$(w_1,w_2,w_3)^T=A_{C_2}^3 \cdot (w_1,w_2,w_3)^T$$ 
Let 
$$
v_t:=t_1 v_1 + t_2 v_2 + t_3 v_3
$$ 
where ($v_1,v_2,v_3$ - eigenvectors from proposition \ref{C2eigen}). It suffices to show that for every $\epsilon_X,\epsilon_Y$ from the yellow region, every compact interval $I_2,I_3 \subset\mathbb R$ and $t_2 \in I_2, t_3 \in I_3$, there exist $T_{I_2,I_3}>0$, big enough, such that for every $t \in (T_{I_2,I_3}, + \infty)$,  we have convergence to $+\infty$ on all of the coordinates of the vectors $A_{C_2}^n(v_t)$ and $\vartheta_{26}(A_{C_2}^n(v_t))$ (as well for $n=1$ their coordinates have to be big enough). Now 
$$ A_{C_2}^n(v_t)= \lambda_1^n t_1 v_1 +\lambda_2^n t_2 v_2 +\lambda_3^n t_3 v_3
$$ 
and in logarithmic coordinates 
$$
\vartheta_{26}(A_{C_2}^n(v_t)) = \vartheta_{26}(\lambda_1^n t_1 v_1 +\lambda_2^n t_2 v_2 +\lambda_3^n t_3 v_3)= 
$$
$$
=\lambda_1^n t_1 (v_{13}, v_{11}+\frac{1-\epsilon_X}{2} v_{13}, v_{12}- \frac{1+\epsilon_Y}{2} v_{13})^T +\lambda_2^n t_2 \vartheta_{26}(v_2) +\lambda_3^n t_3 \vartheta_{26}(v_3)
$$
Since for $\epsilon_X,\epsilon_Y$ from the yellow region the following inequalities hold: 
$$
v_{11}, \ v_{12}, \ v_{13}, \ v_{12}- \frac{1+\epsilon_Y}{2} v_{13}, \ v_{11}+\frac{1-\epsilon_X}{2} v_{13}>0
$$ 
and
$$
\lambda_1 > max(1, |\lambda_2|, |\lambda_3|)
$$
then while passing to the limit with $n \rightarrow +\infty$, the leading terms are those multiplied by $\lambda_1^n t_1$. 
\newline
Hence, cycle $C_2$ is not almost completely unstable, however it cannot be essentially asymptotically stable since the domain of every second local map associated to $C_2$ cycle is a thin cusp, which means that the signs of the local indices alternate through the $C_2$ cycle. 
\end{proof}
\end{thm}

\begin{rmq}
In the same way we can prove about cycle $C_0$ weaker statement than the one from theorem \ref{C0eas}. 
\end{rmq}


\subsection{Stability analysis of the cycles}
In this subsection, we describe bifurcations happening in the sytem with $\epsilon_X, \epsilon_Y$ varying. In the coordinates $u_i:=- \ln(z_i) \in \mathbb R_+$, $i=1,2,3$ and $z_i \in (0,1)$, let us denote 
$$H:=\{ u \ | \ u_1 +u_2 +u_3 =1 \}, \ \ H_+ := H \cap \mathbb R_+^3 
$$
and the projection $PH$ onto H, that is 
$$
PH(u):= \frac{u}{u_1+u_2+u_3}
$$
wherever defined.
\\ \\
Following propositions are just consequences of the prop. \ref{C0eigen}, \ref{C2eigen} and proof of the proposition \ref{C2attract}. 

\begin{prop}
Projection dynamics on $H$ of the cycle $C_2$: 
\begin{enumerate}
\item -for $\epsilon_X> \epsilon_Y$ in the  yellow region and every $u \in \mathbb R_+^3$: 
$$PH(A_{C_2}^n u) \rightarrow PH(v_1) \ \ \  ( \ A_{C_2}^n u \rightarrow +\infty \  on \ \  coordinates)$$
with $n \rightarrow +\infty$ and $u=t_1 v_1 + t_2 v_2 +t_3 v_3 \in \mathbb R_+^3$ like in proposition \ref{C2attract}.
\item -for $\epsilon_X> \epsilon_Y$ in the white region and $u=t_1 v_1 + t_2 v_2$, with  $n \rightarrow +\infty$: 
$$PH(A_{C_2}^n u) \rightarrow PH(v_1)  \ \ ( \ A_{C_2}^n u \rightarrow +\infty \ \ on \  coordinates)$$ 
-moreover for $u=t_1 v_1 + t_3 v_3$, with $n \rightarrow +\infty$: 
$$PH(A_{C_2}^{-n} u) \rightarrow PH(v_1) \ ( \  A_{C_2}^{-n} u \rightarrow 0 \ \  on \ coordinates)$$ 
and every $u=t_1 v_1 + t_2 v_2 + t_3 v_3$: 
$$A_{C_2}^{-2n+1} u \rightarrow - \infty \ \ and \ \ A_{C_2}^{-2n} u \rightarrow +\infty \ \ on \ \ coordinates$$
\newline
-if $t_3\neq 0$ and $PH(u) \in H_+$, then there exist $N>0$ such that: 
$$PH(A_{C_2}^{2N+1} u) \notin H_+$$ 
\item -for $\epsilon_X> \epsilon_Y$ from the red region and $u=t_1 v_1 + t_3 v_3$:
$$PH(A_{C_2}^{-n} u) \rightarrow PH(v_1) \ \  and  \ \ A_{C_2}^{-n} u \rightarrow +\infty \ \  on \ \ coordinates$$
-for $\epsilon_X> \epsilon_Y$ from the light red region and $u=t_1 v_1 + t_2 v_2 + t_3 v_3$:
$$PH(A_{C_2}^{-n} u) \rightarrow PH(v_1) \ \  and  \ \ A_{C_2}^{-n} u \rightarrow +\infty \ \  on \ \ coordinates$$
\end{enumerate}
\end{prop}

\begin{prop}
Projection dynamics on $H$ of the cycle $C_0$: 
\begin{enumerate}
\item -for $\epsilon_X+ \epsilon_Y <0$: 
$$PH(A_{C_0}^n u) \rightarrow PH(v_1) \ \  and  \ \ A_{C_0}^n \rightarrow +\infty \ \ on \ \  coordinates$$  
\item -for $\epsilon_X+ \epsilon_Y >0$ and $u=t_1 v_1 + t_2 v_2 +t_3 v_3$:
$$PH(A_{C_0}^{-n} u) \rightarrow PH(v_1) \ \ and \ \ A_{C_0}^{-n} u \rightarrow +\infty \ \  on \ \ coordinates$$ 
\end{enumerate}
\end{prop}



From the above propositions we can conclude that:
\begin{prop}
$ \ $
\begin{enumerate}
\item for $\epsilon_X > \epsilon_Y$ in the yellow region:
	$dim(W^s(C_2))=4$ and $dim(W^u(C_0))=4$
\item for $\epsilon_X > \epsilon_Y$ in the white region:
	$dim(W^s(C_2))=3$ and $dim(W^u(C_0))=4$
\item for $\epsilon_X > \epsilon_Y$ in the red region:
	$dim(W^u(C_2))=3$ and $dim(W^s(C_0))=4$
\item for $\epsilon_X > \epsilon_Y$ in the light red region:
	$dim(W^u(C_2))=4$ and $dim(W^s(C_0))=4$
\\ \\
Note that for $\epsilon_X<\epsilon_Y$ there are no stable/unstable manifolds of $C_2$, since $v_{12}<0$.
\end{enumerate}
\end{prop}

\begin{figure}
\centering
\begin{tikzpicture}[xscale=1,>=latex]
\coordinate (A) at (-2,0);
\coordinate (B) at (2,0);
\coordinate (C) at (0,3.46410162);
\coordinate (H1) at (0,1.15470054);
\draw (A)-- (B) -- (C) --cycle;
\draw[->] (A) -- (-3,-0.57735027);
\draw[->] (B) -- (3,-0.57735027);
\draw[->] (C) -- (0,4.3);
\draw[decoration={markings, mark=at position 0.5 with {\arrow{>}}}, postaction={decorate}] (1.2,1) -- (0.5,1);
\draw[decoration={markings, mark=at position 0.5 with {\arrow{>}}}, postaction={decorate}] (-0.2,1) -- (0.5,1);
\draw[decoration={markings, mark=at position 0.5 with {\arrow{>}}}, postaction={decorate}] (0.1,1.6) -- (0.5,1);
\draw[decoration={markings, mark=at position 0.5 with {\arrow{>}}}, postaction={decorate}] (0.8,0.55) -- (0.5,1);
\coordinate (A2) at (4,0);
\coordinate (B2) at (8,0);
\coordinate (C2) at (6,3.46410162);
\draw (A2)-- (B2) -- (C2) --cycle;
\draw[->] (A2) -- (3,-0.57735027);
\draw[->] (B2) -- (9,-0.57735027);
\draw[->] (C2) -- (6,4.3);
\draw[decoration={markings, mark=at position 0.5 with {\arrow{>}}}, postaction={decorate}] (5.7,1.2)--(6.6,1.2);
\draw[decoration={markings, mark=at position 0.5 with {\arrow{>}}}, postaction={decorate}] (5.7,1.2)--(5,1.2);
\draw[decoration={markings, mark=at position 0.5 with {\arrow{>}}}, postaction={decorate}] (5.3,1.8) -- (5.7,1.2);
\draw[decoration={markings, mark=at position 0.5 with {\arrow{>}}}, postaction={decorate}] (6,0.75) -- (5.7,1.2);
\coordinate (A3) at (10,0);
\coordinate (B3) at (14,0);
\coordinate (C3) at (12,3.46410162);
\draw (A3)-- (B3) -- (C3) --cycle;
\draw[->] (A3) -- (9,-0.57735027);
\draw[->] (B3) -- (15,-0.57735027);
\draw[->] (C3) -- (12,4.3);
\draw[decoration={markings, mark=at position 0.5 with {\arrow{>}}}, postaction={decorate}] (12.2,1.2)--(12.8,1.2);
\draw[decoration={markings, mark=at position 0.5 with {\arrow{>}}}, postaction={decorate}] (12.2,1.2)--(11.6,1.2);
\draw (13.3,2.95) -- (12.2,1.2);
\draw (11.2,-0.5) -- (12.2,1.2);
\draw[dashed,decoration={markings, mark=at position 0.5 with {\arrow{>}}}, postaction={decorate}] (11.2,-0.5) to[out=10,in=-75] (13.3,2.95);
\draw[dashed,decoration={markings, mark=at position 0.5 with {\arrow{>}}}, postaction={decorate}] (12.6,1.8) to[out=190,in=150] (11.2,-0.5);
\draw[dashed,decoration={markings, mark=at position 0.5 with {\arrow{>}}}, postaction={decorate}] (11.7,0.321428571429) to[out=30,in=-50] (12.6,1.8);
\draw[dashed,decoration={markings, mark=at position 0.7 with {\arrow{>}}}, postaction={decorate}] (12.4,1.47142857143) to[out=180,in=130] (11.7,0.321428571429);
\node at (0.9,1.3) {\scriptsize $PH(v_1)$};
\node at (6.2,1.4) {\scriptsize $PH(v_1)$};
\node at (11.7,1.8) {\scriptsize $PH(v_1)$};
\end{tikzpicture}
\caption{Projective dynamics for cycle $C_2$ and $\epsilon_X>\epsilon_Y$ from yellow, white and red regions respectively.}
\end{figure}
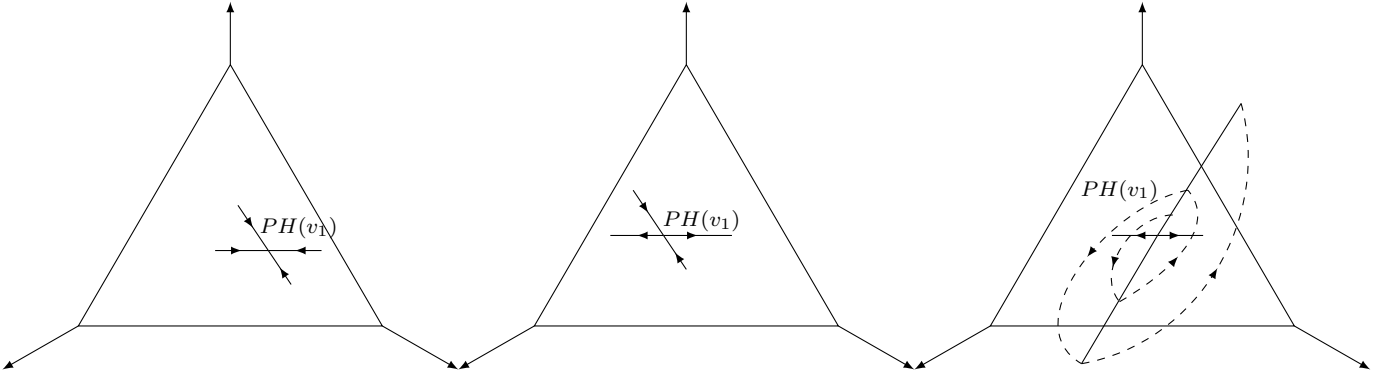

\begin{figure}
\centering												
\begin{tikzpicture}[xscale=1,>=latex]
\coordinate (A) at (-2,0);
\coordinate (B) at (2,0);
\coordinate (C) at (0,3.46410162);
\coordinate (H1) at (0,1.15470054);
\draw (A)-- (B) -- (C) --cycle;
\draw[->] (A) -- (-3,-0.57735027);
\draw[->] (B) -- (3,-0.57735027);
\draw[->] (C) -- (0,4.3);
\draw[decoration={markings, mark=at position 0.5 with {\arrow{>}}}, postaction={decorate}] (1,0.7) -- (0.3,0.7);
\draw[decoration={markings, mark=at position 0.5 with {\arrow{>}}}, postaction={decorate}] (-0.5,0.7) -- (0.3,0.7);
\draw[decoration={markings, mark=at position 0.5 with {\arrow{>}}}, postaction={decorate}] (-0.1,1.3) -- (0.3,0.7);
\draw[decoration={markings, mark=at position 0.5 with {\arrow{>}}}, postaction={decorate}] (0.6,0.25) -- (0.3,0.7);
\coordinate (A2) at (4,0);
\coordinate (B2) at (8,0);
\coordinate (C2) at (6,3.46410162);
\draw (A2)-- (B2) -- (C2) --cycle;
\draw[->] (A2) -- (3,-0.57735027);
\draw[->] (B2) -- (9,-0.57735027);
\draw[->] (C2) -- (6,4.3);
\draw[decoration={markings, mark=at position 0.7 with {\arrow{>}}}, postaction={decorate}] (6.5,1.8)--(6.9,1.8);
\draw[decoration={markings, mark=at position 0.5 with {\arrow{>}}}, postaction={decorate}] (6.5,1.8)--(5.8,1.8);
\draw[decoration={markings, mark=at position 0.5 with {\arrow{>}}}, postaction={decorate}] (6.5,1.8) -- (6.1,2.4);
\draw[decoration={markings, mark=at position 0.5 with {\arrow{>}}}, postaction={decorate}] (6.5,1.8) -- (6.8,1.35);
\node at (0.8,1) {\scriptsize $PH(v_1)$};
\node at (6.2,1.4) {\scriptsize $PH(v_1)$};
\end{tikzpicture}
\caption{Projective dynamics for the cycle $C_0$ and parameters satisfying $\epsilon_X+\epsilon_Y<0$, $\epsilon_X+\epsilon_Y>0$ respectively.}
\end{figure}
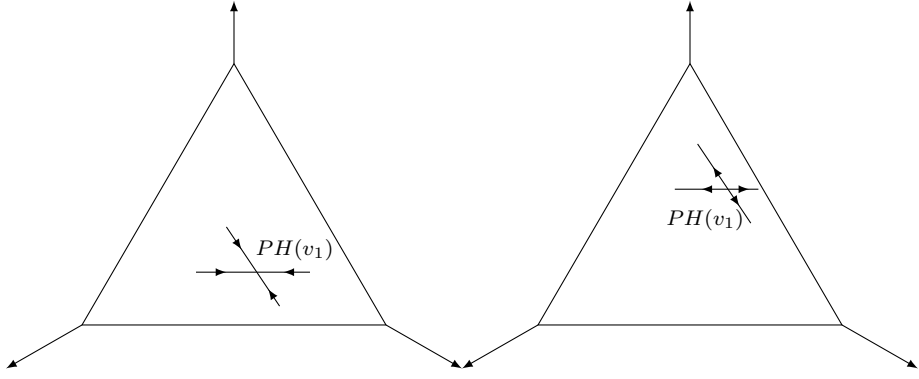

\begin{rmq}
Analogical conclusions follow for $C_1$ cycle after swapping $\epsilon_X$ with $\epsilon_Y$.
\end{rmq}
We end this subsection with a concluding remark.
\begin{rmq}
With parameters $\epsilon_X > \epsilon_Y$ passing from the red region to white region and finally to the yellow region, respectively - we observe a transition from essential asymptotic stability of $C_0$ through the case of Hamiltonian dynamics to the total disappearance of $C_0$ attraction properties and appearance of the 3-dimensional set attracted by $C_2$ cycle, 
then appearance of 4-dimensional set attracted by $C_2$ cycle. 
\newline
When crossing the line $\epsilon_X=\epsilon_Y$ within white region or from the yellow to green region we observe that the set being attracted to $C_2$ cycle disappear and a new one attracted by $C_1$ cycle is born.
\end{rmq}

\section{Shadowing between the ties, lack of infinite switching}\label{sledzenie}
As in the paper by Aguiar\&Castro one can easily prove that there is switching along the heteroclinic orbits in the heteroclinic network. However in this section we point out that their result about infinite switching is incorrect. We provide some examples demonstrating that even attainability of some finite itineraries being infinitely close to the heteroclinic network depends strongly on parameter values $(\epsilon_X, \epsilon_Y)$.
\\ \\
Following Aguiar\&Castro, we introduce the notion of finite switching: 
\begin{df}
We say that there is switching along the heteroclinic orbit - from equilbrium $X$ to equilbrium $Y$, if for all equilbrias, say $Z_1$,...,$Z_k$, for which there is a heteroclinic orbit from $Y$ to $Z_1$,...,$Z_k$, and all equilibria $Q_1,...,Q_m$ from which there is a heteroclinic orbit to $X$, there exist a point $q_{i,j} \in \Sigma^{in}_{Q_i \rightarrow X}$ such that:
$$
\Phi_{Y \rightarrow Z_j} \circ Tran_{X \rightarrow Y} \circ \Psi_{X \rightarrow Y} \circ \Phi_{X \rightarrow Y} (q_{i,j}) \in \Sigma^{out}_{Y \rightarrow Z_j}
$$
for all $i=1,...,m$ and $j=1,...,k$.
\end{df}
\begin{df}
If the process of switching along the heteroclinic orbit can be iterated infinitely many times then we say that there is infinite switching near the heteroclinic network.
\end{df}
\begin{rmq}
In generic case of heteroclinic network where heteroclinic orbits are non-transversal intersections of stable/unstable manifolds, one does not necessarily expect infinite switching to happen arbitrarily close to the heteroclinic network. We refer here to the work by Simon, Delshams and Zgliczy\'nski \cite{[7]}, where they describe possible scenarios of switching along heteroclinic orbits in the non-transversal case and give explanation what one should expect at most from switching and about length of the itinerary, depending e.g. on the dimension of the phase space.
\end{rmq}

\begin{prop}
Finite itineraries following consecutively arbitrarily small neighbourhoods of the points:
\newline
$(P,P) \rightarrow (S,P) \rightarrow (S,S)$ and $(P,P) \rightarrow (S,P) \rightarrow (S,R) \rightarrow (P,R) \rightarrow (P,P)$ 
\newline
are attainable for all $\epsilon_X, \epsilon_Y \in (-1,1)$
\end{prop}
The next proposition presents a counterexample to the Theorem 4.5 about infinite switching claimed in Aguiar\&Castro paper.
\begin{thm}\label{DWLD}
For every $\epsilon_X,\epsilon_Y \in (-1,1)$ the finite itineraries of the form $tie \rightarrow win/loss \rightarrow loss/win \rightarrow tie$ and $tie \rightarrow loss/win \rightarrow win/loss \rightarrow tie$ are not attainable for $h$ - sufficiently small.
\begin{proof}
For an itinerary $(P,P) \rightarrow (S,P) \rightarrow (S,R) \rightarrow (R,R)$ take a point 
$$(z_1,1-h,z_3,z_4) \in \Sigma^{in}_{(P,P)\rightarrow (S,P)} \cap Image( \Psi_{(P,P)\rightarrow (S,P)})
$$ 
Hence $|z_1|<h\cdot (\frac{h}{1-h})^{\frac{3+\epsilon_X}{1-\epsilon_X}}$, $|z_4|<h$. The necessary condition for the point 
$$
Tran_{(S,P) \rightarrow (S,R)} \circ \Psi_{(S,P) \rightarrow (S,R)} \circ \Phi_{(S,P) \rightarrow (S,R)} \circ Tran_{(P,P) \rightarrow (S,P)}(z_1,1-h,z_3,z_4)
$$ 
to lie in the domain of $\Phi_{(S,R) \rightarrow (R,R)}$ is the following inequality: 
$$ 
(\frac{1-h}{h})^{\frac{1+\epsilon_X}{2}} (\frac{z_3}{h})^{\frac{1-\epsilon_X}{2}} |h-z_1| \leq h^{-\frac{1-\epsilon_X}{1+\epsilon_X}} (\frac{h}{1-h})^{\frac{1-\epsilon_X}{1+\epsilon_X}} (\frac{z_1 z_3}{h})^{\frac{2}{1+\epsilon_X}}
$$ 
which is equivalent to
$$
LHS:=\frac{|h-z_1|}{h} \leq h^{-\frac{1-\epsilon_X}{1+\epsilon_X}-1} (\frac{h}{1-h})^{\frac{1-\epsilon_X}{1+\epsilon_X}+\frac{1+\epsilon_X}{2}} (\frac{z_3}{h})^{\frac{2}{1+\epsilon_X} - \frac{1-\epsilon_X}{2}} z_1 ^{\frac{2}{1+\epsilon_X}}=:RHS
$$
Note that 
$$LHS \in \Big( 1-(\frac{h}{1-h})^{\frac{3+\epsilon_X}{1-\epsilon_X}}, 1 \Big)
$$ 
and
$$RHS \in \Big( 0,h^{-\frac{1-\epsilon_X}{1+\epsilon_X}-1+ \frac{2}{1+\epsilon_X}} (\frac{h}{1-h})^{\frac{1-\epsilon_X}{1+\epsilon_X}+\frac{1+\epsilon_X}{2} +\frac{3+\epsilon_X}{1-\epsilon_X} \cdot \frac{2}{1+\epsilon_X}} \Big) =
$$ 
$$
\Big( 0,(\frac{h}{1-h})^{\frac{1-\epsilon_X}{1+\epsilon_X}+\frac{1+\epsilon_X}{2} +\frac{3+\epsilon_X}{1-\epsilon_X} \cdot \frac{2}{1+\epsilon_X}}) \Big)
$$ 
To see that this inequality cannot be satisfied it suffices to say that  
$$(\frac{h}{1-h})^{\frac{1-\epsilon_X}{1+\epsilon_X}+\frac{1+\epsilon_X}{2} +\frac{3+\epsilon_X}{1-\epsilon_X} \cdot \frac{2}{1+\epsilon_X}} < 1-(\frac{h}{1-h})^{\frac{3+\epsilon_X}{1-\epsilon_X}}
$$ 
for $h$-small enough. 
\end{proof}
\end{thm}

\begin{prop}
For every 
$(\epsilon_X,\epsilon_Y) \in (-1,0) \times (-1,1)$, the finite itinerary 
$$
(P,P) \rightarrow (P,S) \rightarrow (R,S) \rightarrow (R,P) \rightarrow (S,P) \rightarrow (S,R) \rightarrow (R,R)
$$ 
is not attainable for sufficiently small $h$.
\begin{proof}
For 
$$
(z_1,z_2,z_3,1-h) \in \Sigma_{(P,P) \rightarrow (P,S)}^{in} \cap \Psi_{(P,P) \rightarrow (P,S)}(\Sigma_{(P,P) \rightarrow (P,S)}^{out})
$$
we have $|z_1|<h$, $|z_2|<h$, $|z_3|<h \ (\frac{h}{1-h})^{\frac{3+\epsilon_Y}{1-\epsilon_Y}}$. 
If 
\begin{multline*}
(w_1,w_2,w_3,w_4)= Tran_{(S,P) \rightarrow (S,R)} \circ \Psi_{(S,P) \rightarrow (S,R)} \circ \Phi_{(S,P) \rightarrow (S,R)} \circ Tran_{(R,P) \rightarrow (S,P)} \circ \Psi_{(R,P) \rightarrow (S,P)} \circ 
\\
\circ \Phi_{(R,P) \rightarrow (S,P)} \circ Tran_{(R,S) \rightarrow (R,P)} \circ \Psi_{(R,S) \rightarrow (R,P)} \circ \Phi_{(R,S) \rightarrow (R,P)} \circ Tran_{(P,S) \rightarrow (R,S)} \circ \Psi_{(P,S) \rightarrow (R,S)} \circ
\\
\circ
\Phi_{(P,S) \rightarrow (R,S)} \circ Tran_{(P,P) \rightarrow (P,S)}(z_1,z_2,z_3,1-h)
\end{multline*}
then the necessary condition for existence of this connection is $|w_2| < h (\frac{|w_1|}{h})^{\frac{2}{1+\epsilon_X}}$, which is equivalent to 
\begin{multline*}
\Big( 1 - (\frac{1-h}{h})^{\epsilon_X} (\frac{h}{|z_1|})^{\frac{1 + \epsilon_X}{2}} \frac{|z_2|}{h} \Big)^{p_1} \cdot \Big(1 - \frac{|z_3|}{h} \Big)^{p_2} < 
\\
< (1 - h)^{\frac{
 \epsilon_X^2 (-5 + 8 \epsilon_Y - 3 \epsilon_Y^2) + \epsilon_X^3 (15 - 8 \epsilon_Y + \epsilon_Y^2) + 
  \epsilon_X (13 - 24 \epsilon_Y + 3 \epsilon_Y^2) - 3 (5 - 8 \epsilon_Y + 3 \epsilon_Y^2)}{
 16 (1 + \epsilon_X)}} \cdot 
h^{-\frac{( 17 + \epsilon_X - \epsilon_X^3 (-3 + \epsilon_Y) + 9 \epsilon_Y - 3 \epsilon_X \epsilon_Y + 3 \epsilon_X^2 (1 + \epsilon_Y)}{8 (1 + \epsilon_X))}} \cdot 
\\
|z_1|^{\frac{ 5 - 3 \epsilon_X (-3 + \epsilon_Y)^2 - \epsilon_X^3 (-3 + \epsilon_Y)^2 + 6 \epsilon_Y + 9 \epsilon_Y^2 + \epsilon_X^2 (-9 + 2 \epsilon_Y + 3 \epsilon_Y^2)}{16 (1 + \epsilon_X)}} \cdot 
|z_2|^{\frac{ 5 + 4 \epsilon_X - \epsilon_X^2 (-3 + \epsilon_Y) - 3 \epsilon_Y}{4 (1 + \epsilon_X)}} |z_3|^{\frac{3 + \epsilon_X^2}{2 + 2 \epsilon_X}} 
\end{multline*}
\newline
All of the powers of $z_1,z_2,z_3$ on the RHS are positive for $(\epsilon_X,\epsilon_Y) \in (-1,0) \times (-1,1)$, and for these parameters $p_1, p_2$ belong to some compact interval in $\mathbb R$. Hence RHS attains maximum for $z_1=h$, $z_2=h$ and $z_3=h (\frac{h}{1-h})^{\frac{3+\epsilon_Y}{1-\epsilon_Y}}$ and is equal to 
$$
\big( \frac{h}{1-h} \big)^{-\frac{87 - 15 \epsilon_Y + 33 \epsilon_Y^2 - 9 \epsilon_Y^3 + 
 \epsilon_X^2 (29 - 5 \epsilon_Y + 11 \epsilon_Y^2 - 3 \epsilon_Y^3) + 
 \epsilon_X^3 (-15 + 23 \epsilon_Y - 9 \epsilon_Y^2 + \epsilon_Y^3) + 
 \epsilon_X (-13 + 37 \epsilon_Y - 27 \epsilon_Y^2 + 3 \epsilon_Y^3)}{16 (1 + \epsilon_X) (-1 + 
   \epsilon_Y)}}
$$
The power is bigger than $0$ for $(\epsilon_X,\epsilon_Y) \in (-1,0) \times (-1,1)$, and hence for $h$ small enough, the RHS is always smaller than LHS which is bounded from $0$.
\end{proof}
\end{prop}
In the same way we investigate the attainability of connections between ties, that is itineraries of the form:
$tie \rightarrow$ (part of $C_0$ cycle) $\rightarrow tie$ 
\\ \\
We would like to end this subsection with small remark about the changes of probability that players will follow some finite itinerary depending on its length.
\begin{prop}
For arbitrarily fixed $\epsilon_X, \epsilon_Y \in (-1,1)$, there exists positive number $p>0$, such that for the finite itinerary $\zeta_1, \zeta_2, \zeta_3$ of the type "$tie \rightarrow  win/loss \ (loss/win) \rightarrow tie$" or "$win/loss \ (loss/win) \rightarrow tie \rightarrow win/loss \ (loss/win)$, the measure of the set of points that follow from 'out cross-section' near $\zeta_1$ to "in cross-section" before $\zeta_3$, through "in cross-section" before $\zeta_2$, is less than $(\frac{h}{1-h})^{p}$ times the measure of the set of points that follow from 'out cross-section' near $\zeta_1$ to "out cross-section" after $\zeta_2$ in the direction to $\zeta_3$. On the other hand measure of the points that follow itinerary of the length 4 is in the limit with $h \rightarrow 0$ approximately the same as the measure of points that follow the same itinerary extended by one equilibrium, in the order "$tie \rightarrow  win/loss \ (loss/win) \rightarrow tie$" or "$win/loss \ (loss/win) \rightarrow tie \rightarrow win/loss \ (loss/win)$.
\end{prop}

\subsection{Chaotic switching}\label{chaos}
In Aguiar\&Castro, the main theorem about RSP game was an incorrect statement about existence of infinite switching for all of the parameter values $\epsilon_X, \epsilon_Y \in (-1,1)$. We have already proven that it is not true, however there is still an open question about the existence of chaos in the system and switching which might be in fact described by subshift of finite type on two symbols (up to topological conjugation, e.g. see \cite{[21]}). 
\newline
We have already mentioned results by Delshams, Simon and Zgliczy\'nski \cite{[7]}, which convince us that there has to exist some phenomenon in the system which forces existence 
of chaotic switching. In the sixth section we will provide numerical investigation in order to support our presumptions and intuition about what stays behind the irregular behaviour in the system (\ref{rownanie}). 
\newline
Typically, examples of infinite switching that appear in the literature are associated to the existence of suspended horseshoes in the neighbourhood of the network (see e.g. \cite{[15]}), in particular if there is an equilibrium in the network with complex eigenvalues - this can lead to the situation similar as in Shilnikov homoclinic bifurcation i.e. spiraling chaos. Although it is worth noting that A.Rodrigues \cite{[14]} has given an example of an attracting homoclinic network exhibiting sensitive dependence on initial conditions, where the existence of infinite switching is not associated with suspended horseshoes. However there is a reason why we believe that there exist chaotic switching topologically conjugated to a subshift of finite type for all parameter values $\epsilon_X, \epsilon_Y$ from the white region:
\begin{prop}(S.Castro, S.van Strien \cite{[5]})
\newline 
For $\epsilon_X, \epsilon_Y$ from the white region satisfying $\epsilon_X=\frac{1+3 \epsilon_Y}{1-\epsilon_Y}$ ($\epsilon_X, \epsilon_Y$ on the $\beta$-curve), system (\ref{rownanie}) has switching type, near the heteroclinic network, corresponding to a subshift of finite type.
\end{prop}

\subsection{Finite switching}
One can hope that it would be possible to determine exactly the form of eventual subshift of finite type describing switching near the system for all parameter values $(\epsilon_X, \epsilon_Y)$ from the white region. In fact, it appears that it is hard to say precisely, for given parameter values $(\epsilon_X, \epsilon_Y)$, which finite itineraries are forbidden for orbits staying arbitrarily close to the network. Below we present some examples of finite itineraries about which we can say only when for sure they are not attainable. 
\newline
Shaded regions in the parameter space indicate parameter values $(\epsilon_X, \epsilon_Y)$ for which the finite itinerary (written under the picture) consisting of numbered cross-sections in the quotient space (see figure 4), 
are not attainable. For simplicity, we denote sequence (3, 4) as $C_0$ and (2,1) as $C_2$.

\begin{figure}[!htb]
\minipage{0.20\textwidth}
  \includegraphics[width=\linewidth]{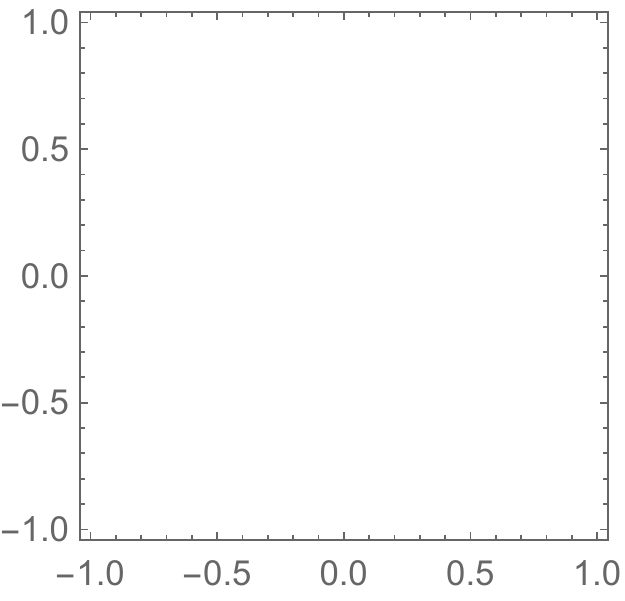}
	{\small (a) 3, 4, 3, 2, 1 
		\newline
		}
\endminipage\hfill
\minipage{0.20\textwidth}
  \includegraphics[width=\linewidth]{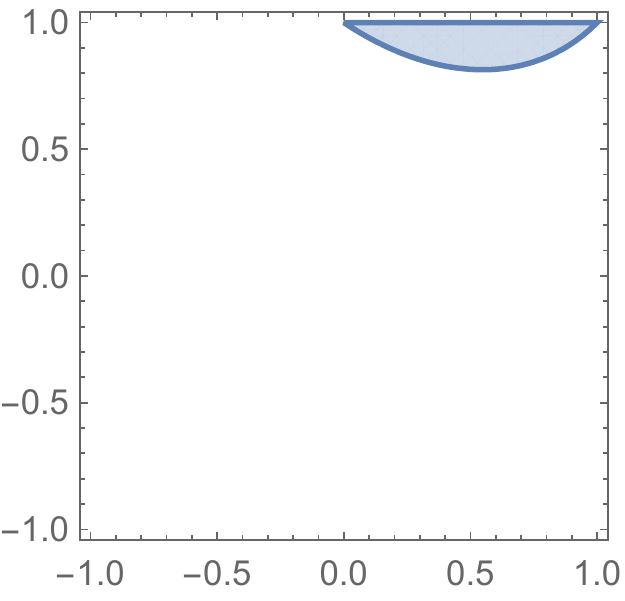}
	{\small (b) 3, 4, 3, 2, 1, 2, 1 
	\newline 
	}
\endminipage\hfill
\minipage{0.20\textwidth}%
  \includegraphics[width=\linewidth]{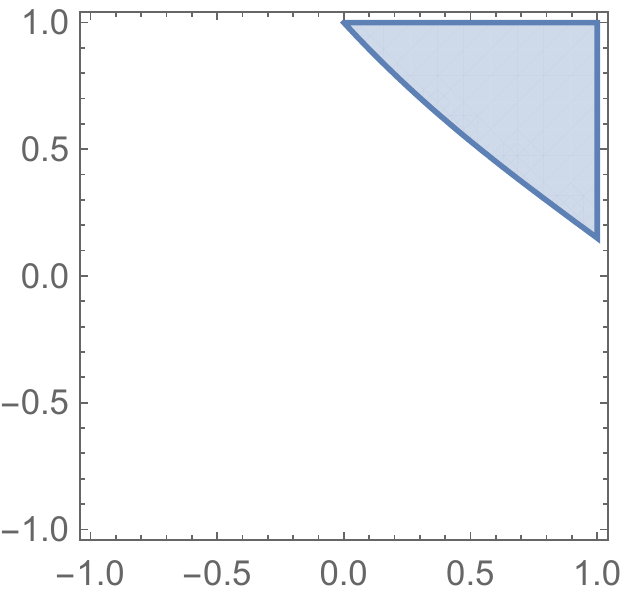}
	{\small (c) 3, 4, 3, 2, 1, 2, 1, 4, 3, 4, 3  }
\endminipage\hfill
\newline 
\minipage{0.20\textwidth}%
  \includegraphics[width=\linewidth]{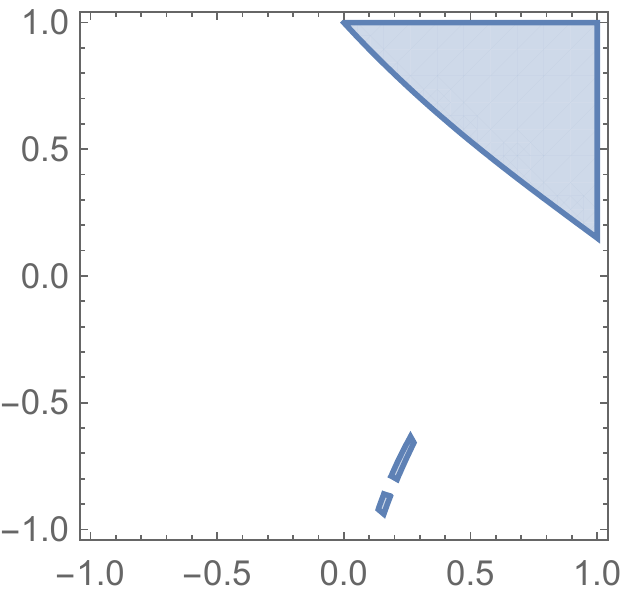}
	{\small (d) 3, 4, 3, 2, 1, 2, 1, 4, 3, 4, 3, 2, 1 }
\endminipage\hfill
\minipage{0.20\textwidth}
  \includegraphics[width=\linewidth]{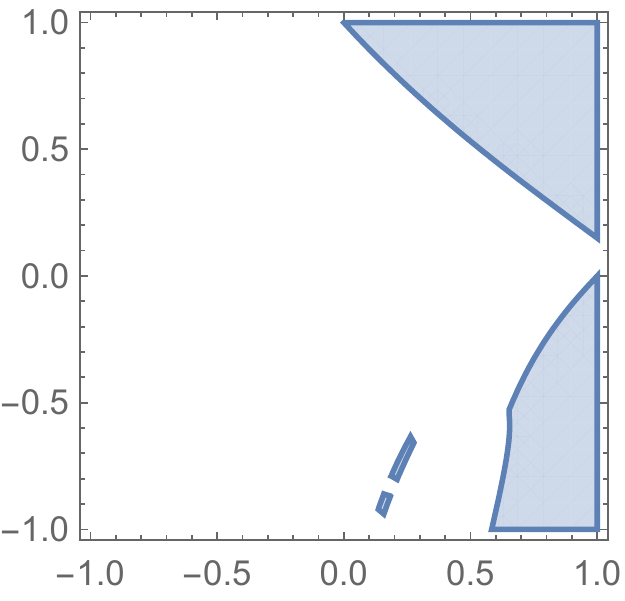}
	{\small (e) 3, 4, 3, 2, 1, 2, 1, 4, 3, 4, 3, 2, 1, 2, 1  }
\endminipage\hfill
\minipage{0.20\textwidth}
  \includegraphics[width=\linewidth]{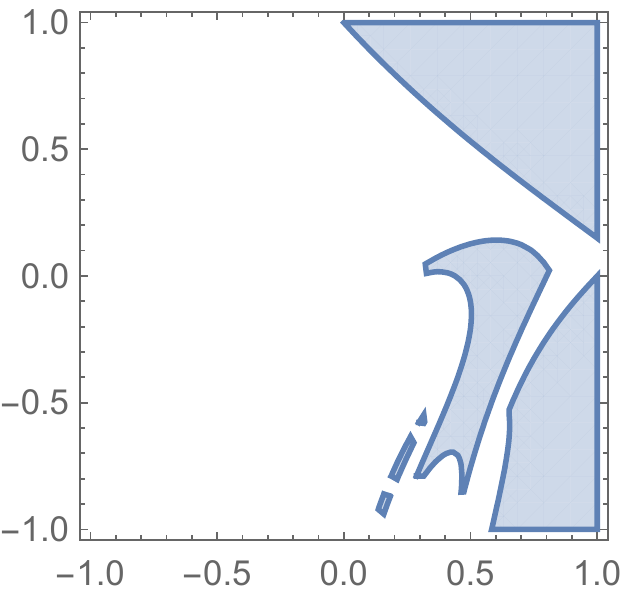}
	{\small (f) 3, 4, 3, 2, 1, 2, 1, 4, 3, 4, 3, 2, 1, 2, 1, 2, 1}
\endminipage\hfill

\end{figure}

\begin{figure}[!htb]
\minipage{0.20\textwidth}
  \includegraphics[width=\linewidth]{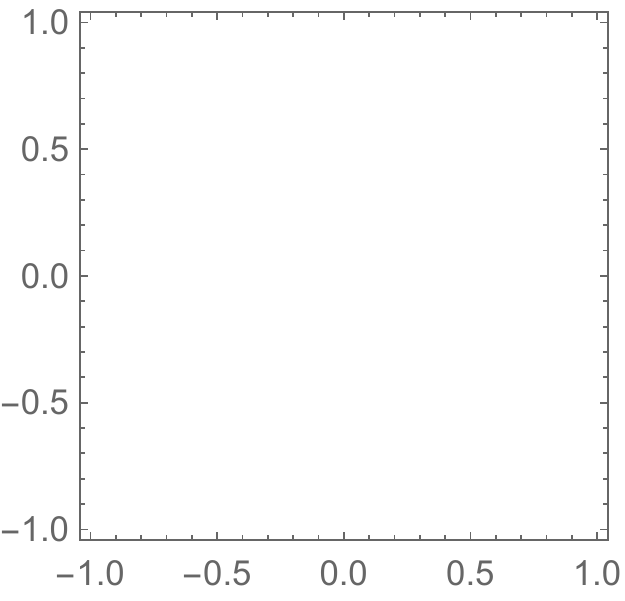}
	{\tiny (a) $C_0$, $C_0$ 
	\newline
	}
\endminipage\hfill
\minipage{0.20\textwidth}
  \includegraphics[width=\linewidth]{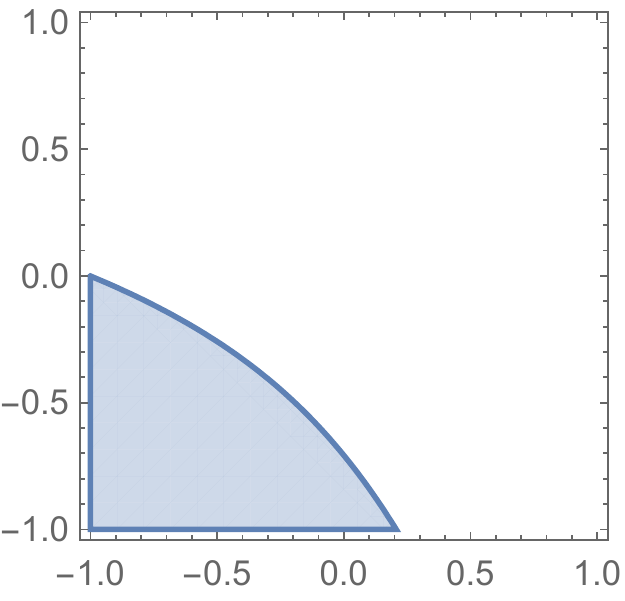}
	{\tiny (b) $C_0$, $C_0$, 3, $C_2$ 
	\newline
	}
\endminipage\hfill
\minipage{0.20\textwidth}%
  \includegraphics[width=\linewidth]{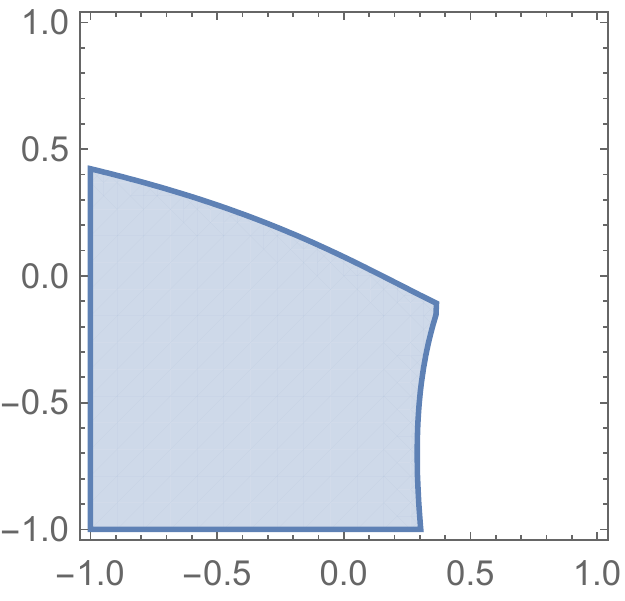}
	{\tiny (c) $C_0$, $C_0$, $C_0$, 3, $C_2$
	\newline 
	}
\endminipage\hfill
\minipage{0.20\textwidth}%
  \includegraphics[width=\linewidth]{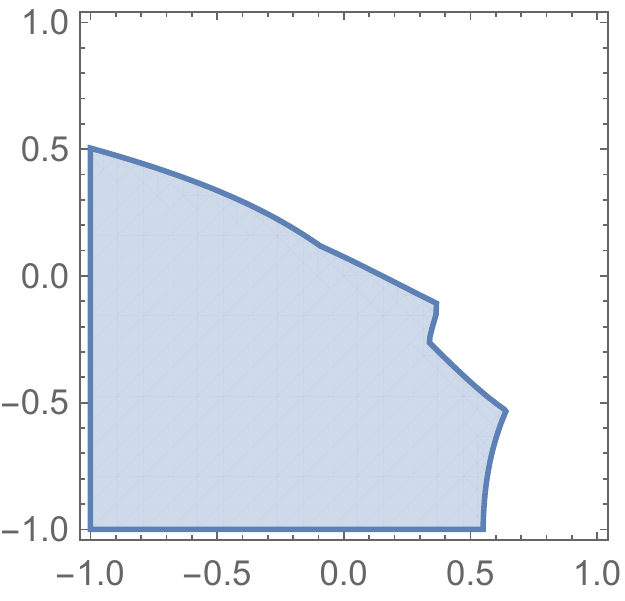}
	{\tiny (d) $C_0$, $C_0$, $C_0$, $C_0$, $C_0$, 3, $C_2$ 
	\newline
	}
\endminipage\hfill
\newline
\minipage{0.20\textwidth}
  \includegraphics[width=\linewidth]{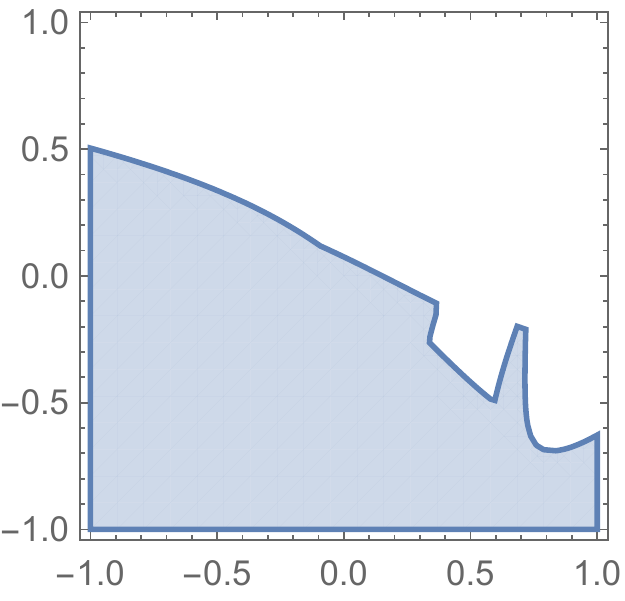}
	{\tiny (e) $C_0$, $C_0$, $C_0$, $C_0$, $C_0$, $C_0$, 3, $C_2$, $C_2$, 4, $C_0$, 3}
\endminipage\hfill
\minipage{0.20\textwidth}
  \includegraphics[width=\linewidth]{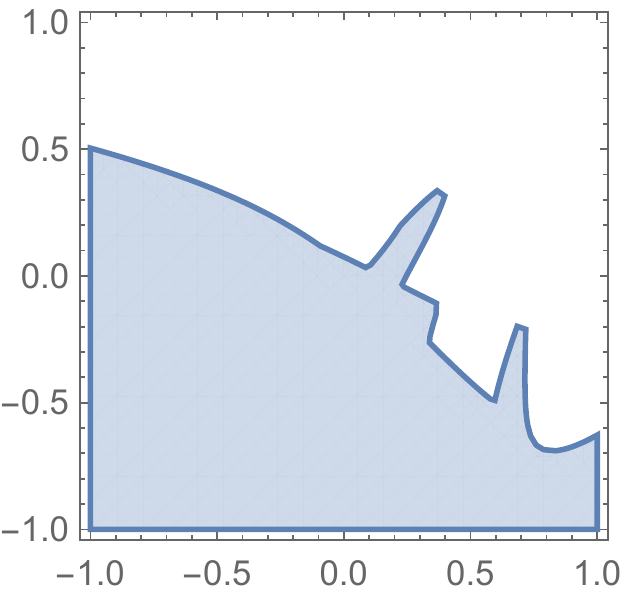}
	{\tiny (f) $C_0$, $C_0$, $C_0$, $C_0$, $C_0$, $C_0$, 3, $C_2$, $C_2$, 4, $C_0$, 3, $C_2$}
\endminipage\hfill
\minipage{0.20\textwidth}
  \includegraphics[width=\linewidth]{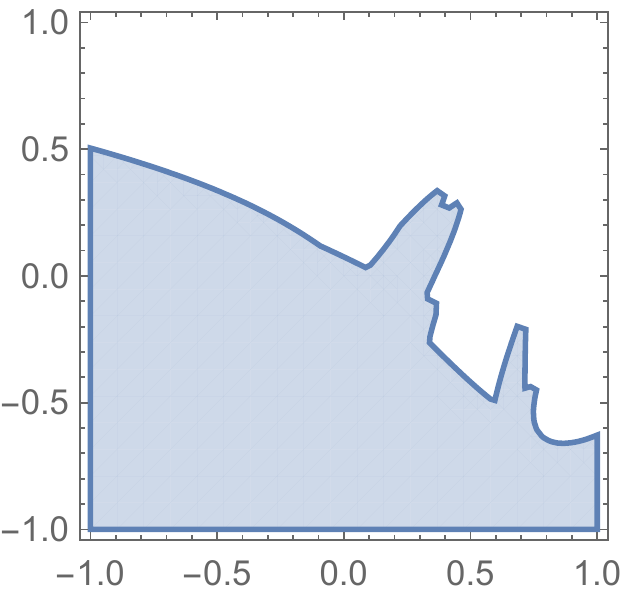}
	{\tiny (f) $C_0$, $C_0$, $C_0$, $C_0$, $C_0$, $C_0$, 3, $C_2$, $C_2$, 4, $C_0$, 3, $C_2$, $C_2$, $C_2$}
\endminipage\hfill

\caption{In these diagrams the parameter space $(\epsilon_X, \epsilon_Y) \in (-1,1)^2$ and various itineraries are included. See main text from section 5.2, for the description.}
\end{figure}

\section{Numerical investigation}
\subsection{Sato's et al. simulations}
In this subsection we remind some of the numerical observations done by Y. Sato et al. (see \cite{[16]}, \cite{[17]}, \cite{[18]}) and in the following subsection we compare it to our findings.  
\begin{enumerate}
\item  With $\epsilon_X=\epsilon_Y=0$ only quasiperiodic tori is observed. For some initial conditions, tori are knotted and form a trefoil.
\item With $\epsilon_X=-\epsilon_Y>0$ Hamiltonian chaos can occur with positive-negative pairs of Lyapunov exponents. The dynamics is very rich in the sense that there are infinitely many distinct behaviors near the Nash equilibrium and a periodic orbit arbitrarily close to any
chaotic one. (see \cite{[6]})
\item When the game is not zero-sum ($\epsilon_X + \epsilon_Y \neq 0$), transients to heteroclinic cycles are observed. On the one hand, there are intermittent behaviors in which the time spent near pure strategies increases subexponentially with $\epsilon_X + \epsilon_Y < 0$ and on the other hand, with $\epsilon_X +\epsilon_Y >0$, chaotic transients persist.
\item When $\epsilon_X + \epsilon_Y < 0$, the behavior is intermittent and orbits are guided by the flow on edges of $\Delta$, which describes a network of possible heteroclinic cycles. Since action ties are not rewarded there is only one such cycle - $C_0$. During the cycle each agent switches between almost deterministic actions in the order $R \rightarrow S \rightarrow P$. The agents are out
of phase with respect to each other and they alternate winning each turn. 
\item When $\epsilon_X + \epsilon_Y > 0$, however, the orbit is an infinitely persistent chaotic transient. Since, in this case, agent $X$ can choose a tie, the cycles are not closed. For example, with $\epsilon_X > 0$, at $(R, P)$, agent $X$ has the option of moving to $(P, P)$ instead of $(S, P)$ with a positive probability. This embeds an instability along the heteroclinic cycle and so
orbits are chaotic.
\item In \cite{[17]}, Lyapunov spectra, for different parameter values satisfying $\epsilon_X + \epsilon_Y=0$ and different initial conditions, were investigated. In conclusion, as an evidence for Hamiltonian chaos for $(\epsilon_X,\epsilon_Y)=(0.5,-0.5)$ was presented the Lyapunov spectrum with the biggest Lyapunov exponent of order $4\cdot 10^{-2}$, and corresponding the second and third were of order $6\cdot 10^{-4}$.  
\end{enumerate}

\subsection{Our numerical observations}\label{wykladniki}

Let us change the coordinates to the ones introduced in \cite{[18]} (eq. (9)), which are useful since being more numerically stable. If $x_j,y_k \neq 0$, we substitute: 
\begin{align}
u_i:= \ln \Big( \frac{x_{i+1}}{x_1} \Big) \ \ \ 
v_i:= \ln \Big( \frac{y_{i+1}}{y_1} \Big)
\end{align}
for $i=1,2$. Then the equation (\ref{rownanie}) reduces to:
\begin{equation}\label{rownanieSN}
\begin{cases}
\dot{u_1}=\frac{-(1+\epsilon_X)- (1-\epsilon_X) e^{v_1}+2 e^{v_2}}{1+e^{v_1}+e^{v_2}} 
\\ 
\dot{u_2}=\frac{-2 +(1-\epsilon_X) e^{v_1} +(1+\epsilon_X) e^{v_2}}{1+e^{v_1}+e^{v_2}}
\\ 
\dot{v_1}=\frac{-(1+\epsilon_Y)- (1-\epsilon_Y) e^{u_1}+2 e^{u_2}}{1+e^{u_1}+e^{u_2}}
\\
\dot{v_2}=\frac{-2 +(1-\epsilon_Y) e^{u_1} +(1+\epsilon_Y) e^{u_2}}{1+e^{u_1}+e^{u_2}}
\end{cases}
\end{equation}
The reverse transformation is given by:
$$
x_1=\frac{1}{1+e^{u_1}+e^{u_2}} \ \ \ 
x_2=\frac{e^{u_1}}{1+e^{u_1}+e^{u_2}} \ \ \
y_1=\frac{1}{1+e^{v_1}+e^{v_2}} \ \ \
y_2=\frac{e^{v_1}}{1+e^{v_1}+e^{v_2}}
$$
Sato's et al. observations were based on simulations done by drawing points randomly from the phase space $\Delta_X \times \Delta_Y$ and integrating their orbits using either low order Runge-Kutta method or some symplectic integrator. However, there are known examples of systems for which usage of symplectic integrators leads to totally different simulations than those obtained by Taylor methods. The following example (Euler's eq. 1760') was provided by Prof. D.Wilczak: 
$$
\begin{cases}
x''=\frac{-2 x}{(\sqrt{x^2 + y^2})^3} - \frac{x-1}{(\sqrt{(x-1)^2 +y^2})^3} 
\\ 
y''=\frac{-2 y}{(\sqrt{x^2 + y^2})^3} - \frac{y}{(\sqrt{(x-1)^2 +y^2})^3}
\end{cases}
$$
One can compare simulated trajectory for the initial point $(x(0),y(0),x'(0),y'(0))= (1.47, 0.8, -0.81, 0)$, and see the differences between simulations using Taylor 30th order method and the symplectic Euler method. One of the reasons for the different portraits is the main disadvantage of the symplectic methods, that is the rounding of errors. For our purposes, we needed to utilise the integrator implemented in CAPD library \cite{[4]}, based on the Taylor method (order of the method can be choosen by the user, we have run the simulations with order from 10 to 30), with the time step fixed $=0.01$ or adjusted after each iteration, in order to minimise the error.
\newline

\subsection{Spectra of Lyapunov exponents}\label{wykladniki}
In order to determine the Lyapunov exponents $\lambda_1 > \lambda_2 > \lambda_3 > \lambda_4$, we utilised the algorithm of Ruelle and Eckmann based on the QR-decomposition method. As the Lyapunov exponents change slightly (difference for the consecutive time steps is of order $10^{-6}$) for the integration time longer than $t_0=2000.0$, we decided to estimate them at time $t_0$.  
\begin{enumerate}
\item For the Hamiltonian case, i.e. $\epsilon_X + \epsilon_Y = 0$, we observe that the biggest Lyapunov exponent $\lambda_1 >0$ is of order $10^{-4}$. As expected for Hamiltonian systems $\lambda_2 \approx -\lambda_3$ are small (of order $10^{-6}$) comparing to $\lambda_4= -\lambda_1$. These results are not compatible with Sato's et al. numerics as we could not find for $\epsilon_X =- \epsilon_Y = 0.5$ any point with the biggest Lyapunov exponent of order $10^{-2}$. In our simulations it was of order $10^{-2}$, in most cases, for the integration time from $t=200.0$ until $t=400.0$. However, for these times, the conditions $\lambda_2 \approx -\lambda_3$, $\lambda_4= -\lambda_1$ were not satisfied even with the precision $10^{-4}$. 
\item For parameter values $\epsilon_X > \epsilon_Y$ from the red and light red region we can observe two different patterns exhibited by the system for various initial conditions:
\newline
- either $\lambda_1>0, \ 0 > \lambda_3, \ \lambda_4$ are of order $10^{-4}$ and $\lambda_2 \approx 0$ is of order $10^{-5}$ or $10^{-6}$, small comparing to $\lambda_1, |\lambda_3|, |\lambda_4|$. For the initial point $(u_1,u_2,v_1,v_2)=(-55.6888,36.3836,-49.0596,-26.4498)$ and parameter values $(\epsilon_X,\epsilon_Y)=(-0.8,-0.9)$, the Lyapunov spectrum is $(2.94135\cdot 10^{-4},2.2451\cdot 10^{-5},-1.4441\cdot 10^{-4},-1.99052\cdot 10^{-4})$.
\newline
- or  $\lambda_1>0$, $\lambda_4 \approx -\lambda_1$ are of order $10^{-4}$ and $\lambda_2, \lambda_3 \approx 0$ are of order $10^{-6}$, $10^{-5}$, respectively (small comparing to $\lambda_1,|\lambda_4|)$. For the initial point $(u_1,u_2,v_1,v_2)=(-2.93175,-71.3114,-53.5168,-87.164)$ and parameter values $(\epsilon_X,\epsilon_Y)=(-0.8,-0.9)$, the Lyapunov spectrum is $(2.10125\cdot 10^{-4},-3.8856\cdot 10^{-6},-1.86203\cdot 10^{-5},-2.114478\cdot 10^{-4})$. Two of the Lyapunov exponents being equal to zero, might indicate that some quantity is preserved by the system. 
\item For parameter values $\epsilon_X > \epsilon_Y$ from the white and yellow regions we can observe three patterns exhibited by the system for various initial conditions. Two of them are identical with those described above, another one is significantly different:
\newline
- $\lambda_1 >0$, $\lambda_4 <0$ are of order $10^{-3}$, $\lambda_2>0$ of order $10^{-4}$ and $\lambda_3 \approx 0$ of order $10^{-5}$ or $10^{-6}$, which is small comparing to $\lambda_1, \lambda_2, |\lambda_4|$. In the white region, for the initial point $(u_1,u_2,v_1,v_2)=(6.755,-27.374,58.7936,-12.1726)$ and parameter values $(\epsilon_X,\epsilon_Y)=(0.7,-0.6)$, the Lyapunov spectrum is $(1.69506\cdot 10^{-3},2.13178\cdot 10^{-4},1.20503\cdot 10^{-5},-1.94714\cdot 10^{-3})$. In the yellow region, for the initial point $(u_1,u_2,v_1,v_2)=(1.19127,-44.7834,99.1464,75.3413)$ and parameter values $(\epsilon_X,\epsilon_Y)=(0.9,0.4)$, the Lyapunov spectrum is $(1.60342\cdot 10^{-3},1.19112\cdot 10^{-4},-6.95594\cdot 10^{-6},-1.74244\cdot 10^{-3})$.  
\end{enumerate}
In each case, we observe that $\lambda_1 + \lambda_2 + \lambda_3 + \lambda_4 < 0$, so the system is dissipative, for all of the parameter values $(\epsilon_X,\epsilon_Y) \in (-1,1)^2$.
\begin{rmq}
At this place we want to thank Yuzuru Sato for his suggestion, that for the parameters from the white region, the biggest Lyapunov exponent is positive and the system (\ref{rownanie}) numerically appears to be chaotic.
\end{rmq}
\subsection{Invariant manifolds of the Nash equilibrium}\label{numerki}
\begin{rmq} 
Note that transformation of the form $(t, \epsilon_X, \epsilon_Y) \ \longmapsto (-t, - \epsilon_X, - \epsilon_Y)$ of the system (\ref{rownanie}), reverses the time and ordering of the cycles $C_0$, $C_1$ and $C_2$. As well it swaps $C_1$ and $C_2$ cycles. Due to this fact, it suffices to conduct forthcoming analysis only for the parameter values $(\epsilon_X, \epsilon_Y)$ from the yellow region and white region below the diagonal (i.e. $\epsilon_X > \epsilon_Y$). 
\end{rmq}
Let us denote $L_1,L_2,L_3,L_4$ - normalised eigenvectors corresponding to the eigenvalues $(\lambda_1, \lambda_2, \lambda_3, \lambda_4)$, from proposition (\ref{NashEq}), for the Nash equilibrium of the system (\ref{rownanieZredukowane}). For $\epsilon_X + \epsilon_Y \neq 0$, vectors having real coefficients: 
\begin{align}
w_1:=L_1 +L_3, \ \ w_2:= i (L_1 -L_3), \ \ w_3:=L_2 +L_4, \ \ w_4:= i (L_2 -L_4)
\end{align}
span the tangent spaces to stable and unstable manifolds, i.e. $span(w_1,w_3)$ and $span(w_2,w_4)$ respectively.
\newline
For the numeric simulations, in order to observe some irregular behaviour or convergence of orbits either to the cycles $C_1$, $C_2$, or $C_0$, we set Poincare sections $S_1:= \{ x_2+y_2=\frac{e^{u_1}}{1+e^{u_1}+e^{u_2}} +\frac{e^{v_1}}{1+e^{v_1}+e^{v_2}} =1 \}$, $S_2:= \{ x_1+y_1=\frac{1}{1+e^{u_1}+e^{u_2}} +\frac{1}{1+e^{v_1}+e^{v_2}} =1 \}$ and we:
\newline
- iterate backward the points from the stable tangent eigenspace $p:=(\frac{1}{3},\frac{1}{3},\frac{1}{3},\frac{1}{3}) + \alpha_1 w_1 + \alpha_3 w_3$ (in $(x_1,x_2,y_1,y_2)$ coordinates) with $\alpha_1, \alpha_3 \in (-\frac{1}{1000},\frac{1}{1000})$, via the Poincare map associated to the sections  $S_1$, $S_2$ 
\newline 
- iterate forward the points from the unstable tangent eigenspace $q:=(\frac{1}{3},\frac{1}{3},\frac{1}{3},\frac{1}{3}) + \alpha_2 w_2 + \alpha_4 w_4$ (in $(x_1,x_2,y_1,y_2)$ coordinates) with $\alpha_2, \alpha_4 \in (-\frac{1}{1000},\frac{1}{1000})$, via the Poincare map associated to the sections  $S_1$, $S_2$ 
\\ \\
Observations on forward integration of points $q$:
\begin{enumerate}
\item for $\epsilon_X$, $\epsilon_Y$ from the white region, 
\newline
- if $\epsilon_X + \epsilon_Y \approx 0$, we observe irregular behaviour, we have not observed any symptoms of convergence within the first 200 iterations of the Poincare map
\newline
- if $\epsilon_X + \epsilon_Y$ is small but not $\approx 0$, the motion is irregular, however we can also observe iterations of $q$ to follow $C_2$ cycle for some time or for very few iterations come closer to $C_0$ cycle, but then come back to the irregularity
\newline
- for $\epsilon_X + \epsilon_Y$ away from $0$, we can observe either irregular behaviour for the whole simulation or it might happen that there are some irregularities at the beginning of iteration process, but the long term behaviour is the convergence to the $C_2$ cycle, this case requires additional investigation as we believe, the behaviour can depend on the sign of $\epsilon_Y$ as well
\newline
- for $(\epsilon_X, \epsilon_Y)$ being close to the boundary with yellow region, we observe fast convergence to the $C_2$ cycle 
\item for $\epsilon_X$, $\epsilon_Y$ from the yellow region, 
\newline
- we do not observe any symptoms of irregular behaviour of the iterations of points $q$, almost all of them converge fastly to the $C_2$ cycle
\end{enumerate}
Observations on backward integration of points $p$ are obviously similar in few aspects:
\begin{enumerate}
\item for $\epsilon_X$, $\epsilon_Y$ from the white region, 
\newline
- if $\epsilon_X + \epsilon_Y \approx 0$, we observe irregular behaviour, we have not observed any symptoms of convergence of the iterations 
\newline
- for $\epsilon_X + \epsilon_Y$ away from $0$, we can observe either irregular behaviour or it might happen that there are some irregularities (especially at the very beginning of iteration process), but the long term behaviour is the convergence to the $C_0$ cycle, this case as well requires additional investigation 
\item for $\epsilon_X$, $\epsilon_Y$ from the yellow region, 
\newline
- we do not observe any symptoms of irregular behaviour of the iterations of points $p$, almost all of them converge fastly to the $C_0$ cycle
\end{enumerate}

We do not precise what do we mean by "fast convergence", $\epsilon_X + \epsilon_Y \approx 0$ or $\epsilon_X + \epsilon_Y$ far away from $0$, as this is beyond our goals and scope of this paper. 
The purpose of conducting these numerical simulations was to explain: what kind of phenomenon can cause chaotic behaviour in the system or even near the heteroclinic network (switching possibly)? After introduction to this question in the subsection (\ref{chaos}), we can now pose our questions and conjectures. 

\section{Conclusions and open questions}
For $(\epsilon_X, \epsilon_Y)$ from the yellow (resp. red) region away from antidiagonal almost all of the points from the stable (unstable) manifold after backward (forward) iteration by the Poincare map converge fastly to the cycles $C_0$ or $C_2$ what would suggest that there are not only heteroclinic connections Nash equilibrium -- cycles, but that the intersections of stable (unstable) manifolds of the cycles with unstable (stable) manifolds of the Nash equilibrium are $2$-dimensional. Do these intersections persist (presumably losing 1 dimension) when parameter values vary, especially when they are in white region?
\newline 
Although essential asymptotic stability of $C_0$ cycle for $\epsilon_X + \epsilon_Y <0$ does not exclude possibility of irregular behaviour in the system, it means however that closely to the equilibria $loss/win$, $win/loss$ there is not much free space for switching from $C_0$ to the other cycles. Hence it seems reasonable, why Sato et al. did not observe numerically irregular behaviour, but mentioned in their papers only orbits guided by $C_0$ cycle.
\newline
The situation is different, however, for the cycle $C_2$ and parameter values in yellow region, since basins of attraction at every second equilibrium in this cycle are thin cusps, hence there is a lot of space to switch between the cycles. This, together with the fact that stable manifold of the cycle $C_2$ is $3$-dimensional for $\epsilon_X$, $\epsilon_Y$ in the white region, would explain why Sato et al. concluded that for $\epsilon_X + \epsilon_Y > 0$ observed orbits are chaotic and the cycles are not closed, even for the $\epsilon_X$, $\epsilon_Y$ in the yellow region.
\newline
We have to point out that propositions \ref{C0eigen}, \ref{C2eigen}, guarantee locally strong contraction of the attracted open sets of points towards cycles $C_0$ and $C_2$, for $\epsilon_X + \epsilon_Y$ being very small or very big, respectively. Hence it seems plausible that the iterations of points from the tangent eigenspaces of the Nash equilibrium are converging fastly to the cycles - where we would look for the forward and backward heteroclinics from the Nash equilibrium. On the other hand, for $\epsilon_X + \epsilon_Y$ close to $0$, this local contraction is weak (for the $C_2$ cycle - of the 3-dimensional set of points), so our numerical observations about lack of convergence of points $p$ and $q$ towards the cycles, do not stay in the contradiction with e.g. essential asymptotic stability of the $C_0$ cycle. Irregular behaviour which we observe for these parameter values might come from the fact that the system then is close to the Hamiltonian one (the case $\epsilon_X + \epsilon_Y =0$), where we observe (and it was broadly described in \cite{[16]}, \cite{[17]}, \cite{[18]}, \cite{[6]}) Hamiltonian chaos. 
\newline
The irregular behaviour we observe for $(\epsilon_X$, $\epsilon_Y)$ from the white region and the possible switching we expect to happen, might be as well a consequence of existence of homoclinic orbit to the cycles $C_0$ and $C_2$ (superhomoclinic orbit \cite{[20]}), which might not pass very close to the Nash equilibrium, so we did not have the chance to notice it. 
\newline
On the other hand, our numerical experiments suggest, that for some parameter values $\epsilon_X$, $\epsilon_Y$, in white region, far away from the anti-diagonal and close to the bifurcation curve between white and yellow region, there exist a forward heteroclinic connection from the neighbourhood of the Nash equilibrium to the cycle $C_2$, and backward heteroclinic to the cycle $C_0$. Unfortunately, the huge obstacle for performing computer assisted proof of such connections (even for one choosen pair of parameter values ($\epsilon_X,\epsilon_Y$)) is the parametrization and analytical estimations of the stable and unstable manifolds of these cycles (see e.g. \cite{[1]}, \cite{[2]}, \cite{[3]}, \cite{[28]}, \cite{[29]}, \cite{[11]}, \cite{[25]}, \cite{[24]}, \cite{[27]}). Existence of these heteroclinics together with possibility of switching from $C_2$ cycle to the $C_0$ at any equilibrium $loss/win$ (so where the first agent loses), would result in chaotic behaviour possibly described by subshift of finite type, i.e. a point starting close to the Nash equilibrium can follow the heteroclinic, goes around $C_2$ cycle for some time, switches and encircles $C_0$ and then comes back to the neighbourhood of Nash equilibrium and repeats. 
\newline 
If there exists such a connection, i.e. superhomoclinic orbit, passing by close to the Nash equilibrium, for the transportation of chaos we would expect it to be a transversal intersection of unstable and stable manifold of $C_0$ and $C_2$ cycles, respectively. Generically it should be 3-dimensional set for the parameter values $(\epsilon_X$, $\epsilon_Y)$ from the white region. Unfortunately this is another obstacle which makes us think that it is very hard to prove its existence or even compute it numerically (see \cite{[10]}).   
\newline
Our belief is that the transversal intersection of unstable and stable manifold of $C_0$ and $C_2$ cycles, respectively, exists for parameter values $(\epsilon_X$, $\epsilon_Y)$ from the white region and with $(\epsilon_X$, $\epsilon_Y)$ crossing the bifurcation curve towards yellow region, it bifurcates to two 2-dimensional surfaces which are intersections of unstable (resp. stable) manifold of cycle $C_0$ ($C_2$) with stable (unstable) manifold of Nash equilibrium. During this bifurcation, 3-dimensional set would disappear, while $C_2$ gains one dimension of stabilty and two 2-dimensional invariant sets would be born.
\newline
It seems that because of the nature of the problem (nonlinear vector field, phase space being a cartesian product of simplices, the dependence of the system on the parameters) the open questions and presumptions stated above, are not easy to tackle analytically nor numerically (when one wants to provide the dependence of the picture for the whole parameter space). Although the description of the heteroclinic network provided in this paper seems to give a good insight into the local behaviour near the network, but in fact it might turn out that it is not very helpful to dispel all of the doubts about chaotic switching, as we know and mentioned in the section 5, the switching depends strongly on the parameter values. Nevertheless, as we strongly believe - from the applications point of view far more important questions concern global behaviour in the system instead of local near the heteroclinic network. That is why, at the very end, we would like to point out that it should not be hard to perform computer assisted proof, for example of existence of chaos (if the chaos really appears), for open set of parameter values $(\epsilon_X$, $\epsilon_Y)$ in the white region (see e.g. \cite{[26]}), especially for the parameters lying on the $\beta$-curve. However this might not give the full picture of the switching happening on the $\beta$-curve and explanation how it is born near the heteroclinic network. Our future work is aimed at this direction. 

\bigskip

\noindent{\bf Acknowledgements.}
I want to thank Prof. Sebastian van Strien, Sofia Castro, Yuzuru Sato, Alexandre Rodrigues, Prof. Piotr Zgliczy\'nski, Prof. Rafael De La Llave and others for the discussions. I would like to express my special gratitude to Prof. Dmitry Turaev for constant supervision while working on this paper, to Prof. Daniel Wilczak for discussions as well as essential help with the CAPD library, and to the Jagiellonian University in Krak\'ow for the hospitality during my visits there.
This research was partially supported by Santander Mobility Grant.


\begin{thebibliography}{99}

\bibitem{[1]}
Baldoma I., Fontich E., De La Llave R., Martin P.:
\newblock  The parameterization method for one-dimensional invariant manifolds of higher dimensional parabolic fixed points.
\newblock {\em Discrete and continuous dynamical systems}, {\bf 17(4)} (2007) 835

\bibitem{[2]}
Cabre X., Fontich E., Llave R. D. L.:
\newblock  The parameterization method for invariant manifolds II: regularity with respect to parameters.
\newblock {\em Indiana University mathematics journal}, {\bf 52(2)} (2003) 329-360

\bibitem{[3]}
Cabre X., Fontich E.,De La Llave R.: 
\newblock  The parameterization method for invariant manifolds III: overview and applications.
\newblock {\em Journal of Differential Equations}, {\bf 218(2)} (2005) 444-515

\bibitem{[4]}
CAPD Library:   
\newblock http://capd.ii.uj.edu.pl/
\newblock {\em } {\bf } 

\bibitem{[5]}
Castro S.B.S.D., van Strien S.: 
\newblock Realising orbits of best response dynamics as time averaged replicator dynamics orbits.
\newblock {\em in preparation}

\bibitem{[6]}
Chawanya T.: 
\newblock A new type of irregular motion in a class of game dynamics systems.
\newblock {\em Prog. Theor. Phys}, {\bf 94} (1995), 163.

\bibitem{[7]}
Delshams A., Simon A., Zgliczy\'nski P.: 
\newblock Shadowing of non-transversal heteroclinic chain.
\newblock {\em in preparation}, {\bf } 

\bibitem{[8]}
Garrido Silva L.: 
\newblock PhD Thesis, University of Porto.
\newblock {\em in preparation}

\bibitem{[9]} 
Hofbauer J.: 
\newblock Evolutionary dynamics for bimatrix games: a hamiltonian system? 
\newblock {\it J. Math. Biol.} {\bf 34} (1996) 675-688

\bibitem{[10]}
Knobloch E., Moore D.R.: 
\newblock Chaotic travelling wave convection. 
\newblock {\em  European J. Mech. B.}, {\bf 10} (1991) 37-42

\bibitem{[11]} 
Kokubu H., Wilczak D., Zgliczy\'nski P.:
\newblock Rigorous verification of cocoon bifurcations in the Michelson system. 
\newblock {\it Nonlinearity 20} {\bf 9} (2007) 2147-2174

\bibitem{[12]}
Lohse A.: 
\newblock Stability of heteroclinic cycles in transverse bifurcations.
\newblock {\em Physica D}  {\bf 310}  (2015),  95-103. 

\bibitem{[13]}
Podvigina O., Ashwin P.: 
\newblock On local attraction properties and a stability index for heteroclinic connections.
\newblock {\em Nonlinearity}, {\bf 24} (2011), 887-929.

\bibitem{[14]}
Rodrigues A.A.P.: 
\newblock Is there switching without suspended horseshoes?
\newblock arXiv:1511.08643v1 

\bibitem{[15]}
Rodrigues A.A.P., Ibanez S.: 
\newblock On the Dynamics Near a Homoclinic Network to a Bifocus: Switching and Horseshoes.
\newblock {\em Int. J. Bif. Chaos}, {\bf 25(11)}, 15300 (19pp), (2015)

\bibitem{[16]}
Sato Y., Akiyama E., Doyne Farmer J.: 
\newblock Chaos in learning a simple two-person game.
\newblock {\em Proc. Natl. Acad. Sci. USA}, {\bf 99} (2002), 4748-4751.
	
\bibitem{[17]}
Sato Y., Akiyama E., Crutchfield J.P.: 
\newblock Stability and diversity in collective adaptation.
\newblock {\em Physica D}, {\bf 210} (2005), 21-57.

\bibitem{[18]}
Sato Y., Crutchfield J.P.: 
\newblock Coupled replicator equations for the dynamics of learning in multiagent systems.
\newblock {\em Physical Review E}, {\bf 67} (2003), 015206(R)

\bibitem{[19]}
Schuster P., Sigmund K., Hofbauer J., Wolff R.: 
\newblock Selfregulation of behavior in animal societies. Part I. Symmetric contests.
\newblock {\em  Biol. Cybern.}, {\bf 40} (1981) 1-8.

\bibitem{[20]}
Shilnikov L.P., Turaev D.: 
\newblock Super-homoclinic orbits and multi-pulse homoclinic loops in Hamiltonian systems with discrete symmetries.
\newblock {\em  Regul. Chaotic Dyn.}, {\bf 2} (1997) 126-138.

\bibitem{[21]}
Srzednicki R., W\'ojcik K.: 
\newblock A geometric method for detecting chaotic dynamics.
\newblock {\em  J. Differential Equations}, {\bf 135} (1997) 66-82.

\bibitem{[22]}
Wilczak D., Zgliczy\'nski P.: 
\newblock Connecting orbits for a singular nonautonomous real Ginzburg-Landau type equation.
\newblock {\em  SIAM Journal on Applied Dynamical Systems Vol. 15}, {\bf 1} (2016) 495-525

\bibitem{[23]}
Wilczak D., Zgliczy\'nski P.: 
\newblock $C^r$-Lohner algorithm.
\newblock {\em  Schedae Informaticae, Vol. 20}, (2011) 9-46

\bibitem{[24]}
Wilczak D., Zgliczy\'nski P.: 
\newblock Heteroclinic Connections between Periodic Orbits in Planar Circular Restricted Three Body Problem - A Computer Assisted Proof.
\newblock {\em  Communications in Mathematical Physics}, 234.1 (2003): 37-75

\bibitem{[25]}
Wilczak D., Zgliczy\'nski P.: 
\newblock Heteroclinic Connections between Periodic Orbits in Planar Circular Restricted Three Body Problem. Part II.
\newblock {\em  Communications in Mathematical Physics}, 259.3 (2005): 561-576.

\bibitem{[26]}
Wilczak D.: 
\newblock Chaos in the Kuramoto-Sivashinsky equations - a computer assisted proof.
\newblock {\em  Journal of Differential Equations, Vol.194}, (2003) 433-459

\bibitem{[27]}
Wilczak D., Serrano S., Barrio R.: 
\newblock Coexistence and dynamical connections between hyperchaos and chaos in the 4D R\"ossler system: a Computer-assisted proof. 
\newblock {\em  SIAM Journal on Applied Dynamical Systems, Vol. 15}, {\bf 1} (2016) 356-390

\bibitem{[28]}
Zgliczy\'nski P., Gidea M.: 
\newblock Covering relations for multidimensional dynamical systems I.
\newblock {\em  J. of Diff. Equations}, {\bf 202} (2004) 32-58

\bibitem{[29]}
Zgliczy\'nski P., Gidea M.: 
\newblock Covering relations for multidimensional dynamical systems II.
\newblock {\em  J. of Diff. Equations}, {\bf 202} (2004) 59-80




\end{thebibliography}
\end{document}